\documentclass[a4paper,11pt,reqno]{amsart}

\usepackage{amsmath,amsfonts,amsthm,color,amssymb, mathrsfs}
\usepackage{xcolor}
\usepackage[T1]{fontenc}
\usepackage{wasysym}
\usepackage{ulem}
\usepackage{stmaryrd} 
\usepackage[left=1in, right=1in, top=1.1in,bottom=1.1in]{geometry}
\usepackage{enumitem}
\usepackage{hyperref}
\hypersetup{
     colorlinks   = true,
     citecolor    = blue,
     linkcolor=blue
}
\usepackage{cancel}
\usepackage{comment}

\usepackage{csquotes}
\usepackage{graphicx}
\usepackage{pstricks}
\usepackage{lmodern}
\numberwithin{equation}{section}

\newcommand{\1}{{\bf 1}}

\newcommand{\yti}{\tilde{y}}

\newcommand{\xiti}{\tilde{\xi}}

\newcommand{\R}{\mathbb R}

\newcommand{\ca}{\mathcal A}
\newcommand{\cb}{\mathcal B}
\newcommand{\cac}{\mathcal C}
\newcommand{\cd}{\mathcal D}
\newcommand{\ce}{\mathcal E}
\newcommand{\cf}{\mathcal F}
\newcommand{\cg}{\mathcal G}
\newcommand{\ch}{\mathcal H}
\newcommand{\ci}{\mathcal I}
\newcommand{\cj}{\mathcal J}
\newcommand{\ck}{\mathcal K}

\newcommand{\bfb}{\mathbf B}

\newcommand{\al}{\alpha}
\newcommand{\der}{\delta}

\newcommand{\ga}{\gamma}
\newcommand{\ka}{\kappa}
\newcommand{\la}{\lambda}

\newcommand{\vp}{\varphi}

\newcommand{\lp}{\left(}
\newcommand{\rp}{\right)}

\def\E{\mathbb E}

\newtheorem{theorem}{Theorem}[section]

\newtheorem{corollary}[theorem]{Corollary}

\newtheorem{definition}[theorem]{Definition}

\newtheorem{lemma}[theorem]{Lemma}
\newtheorem{notation}[theorem]{Notation}

\newtheorem{proposition}[theorem]{Proposition}

\theoremstyle{remark}
\newtheorem{remark}[theorem]{Remark}

\title[Hyperbolic Anderson model]
{Hyperbolic Anderson model 2:\\ Strichartz estimates and Stratonovich setting}

 \author[X. Chen]{Xia Chen}
\address{X. Chen: Department of Mathematics, University of Tennessee,   Knoxville}
\email{xchen@math.utk.edu}

\author[A. Deya]{Aur\'elien Deya}
\address{A. Deya: Institut Elie Cartan, University of Lorraine}
\email{Aurelien.Deya@univ-lorraine.fr}

\author[J. Song]{Jian Song}
\address{J. Song: Research Center for Mathematics and Interdisciplinary Sciences, Shandong University; 
 School of Mathematics, Shandong University}
\email{txjsong@sdu.edu.cn}

\author[S. Tindel]{Samy Tindel}
\address{S. Tindel: Department of Mathematics, 
Purdue University}
\email{stindel@purdue.edu}



\date{}

\begin{document}

\maketitle

%For every $\al\in \R$, we consider the homogeneous Sobolev space $\dot{\ch}^{\al}$ defined by the condition
%$$\|\vp\|_{\dot{\ch}^\al}^2:=\int_{\R^d} d\xi\, |\xi|^{\al} |\cf(\vp)(\xi)|^2 \ < \, \infty\, .$$

\begin{abstract} We study a wave equation in dimension $d\in \{1,2\}$ with a multiplicative space-time Gaussian noise. The existence and uniqueness of the Stratonovich solution is obtained under some conditions imposed on the Gaussian noise. The strategy is to develop some Strichartz type estimates for the wave kernel in weighted Besov spaces, by which we can prove the wellposedness of an associated Young-type  equation. Those Strichartz bounds are of independent interest.
\end{abstract}

\tableofcontents

\section{Introduction}

In \cite{CDST} we have started a long term project aiming at defining wave equations driven by rough noises. More specifically, \cite{CDST} focused on the following Skorohod type equation on $\R_+\times \R^d$ for $d\in \{1,2, 3\}$:
\begin{equation}\label{e:wave-sko}
\frac{\partial^2 u}{\partial t^2} (t,x) =\Delta u(t,x) + u \diamond \dot W(t,x),
\end{equation}
where $\diamond$ stands for the Wick product and where $\dot W$ is a centered Gaussian noise. In \cite{CDST}, the covariance function for $\dot W$ was fractional in time with proper decay in space. That is for $s,t\in\R_{+}$ and $x,y\in\R^{d}$ we had 
\begin{equation}\label{e:cov-sko}
\E[\dot W(s,x) \dot W(t,y)]=|s-t|^{-a_0}\gamma(x-y),
\end{equation}
where $a_0\in[0,1)$ and $\gamma$ is a nonnegative and nonnegative definite (generalized) function admitting a spectral measure $\mu$. Then we obtained that under appropriate conditions on the initial conditions,  the following relation is necessary and sufficient in order to solve~\eqref{e:wave-sko}:
\begin{equation}\label{e:con-sko}
\int_{\R^d} \left(\frac1{1+|\xi|}\right)^{3-a_0}\mu(d\xi) <\infty.
\end{equation}
Notice that \eqref{e:con-sko} quantifies how irregular in space our noise can be according to its regularity in time.  In particular, if the function $\gamma$ in \eqref{e:cov-sko} satisfies the scaling property
\begin{equation}\label{e:gamma-scaling}
\gamma(cx) =  c^{-a}  \gamma(x) 
\end{equation}
for some $a\in(0,d]$, then condition \eqref{e:con-sko} can be recast as 
\begin{equation}\label{e:con-a}
a_0+a<3\,.
\end{equation}

In the current paper, we will contrast the neat condition obtained in \eqref{e:con-sko} with the situation for a pathwise version of the wave equation (interpreted in the Young or Stratonovich sense). Namely, we consider the following Stratonovich type wave equation on $\R_+\times \R^d$  for $d\in\{1,2\}$, 
\begin{equation}\label{e:wave}
\frac{\partial^2 u}{\partial t^2} (t,x) =\Delta u(t,x) + u \dot W(t,x),
\end{equation}
with initial conditions $u(0,x)=u_0(x)$ and $\frac{\partial }{\partial t}u(0,x)=u_1(x)$. The Gaussian noise we consider in \eqref{e:wave} also has the covariance \eqref{e:cov-sko} with $\gamma$ being a nonnegative definite (generalized) function. We prove the following result on the existence and uniqueness of the Stratonovich solution to \eqref{e:wave} (see Proposition \ref{prop:wave-eq} and Corollary \ref{cor:wave-eq} for a precise account).
\begin{theorem}\label{thm:wave-eq}
Assume  $d\in\{1,2\}$ and the spectral measure $\mu$ of $W$ verifies
\begin{equation}\label{e:con-sol'}
\int_{\R^d} \left(\frac{1}{1+|\xi|}\right)^{\rho_d-a_0-\eta}\mu(d\xi)<\infty \, \quad \text{ for some } \eta>0,
\end{equation}
where the parameter $\rho_{d}$ is such that $\rho_1=1$ and $\rho_2=\frac12$. 
 Then, under some  regularity conditions on $u_0$ and $u_1$, there exists a unique Stratonovich solution  to \eqref{e:wave} in a proper path space. 
If  $\gamma$ has the scaling property \eqref{e:gamma-scaling}, condition \eqref{e:con-sol'} is equivalent to \begin{equation}\label{e:con-wave}
\begin{cases}
a_0+a<1, & \text{ if } d=1,\vspace{0.2cm}\\
a_0+a<\frac12, & \text{ if } d=2.
\end{cases}
\end{equation}

\end{theorem}

It is well known that Stratonovich solutions to stochastic PDEs demand a more regular noisy input $\dot W$ than in the Skorohod setting. And indeed in our case, it is readily checked that \eqref{e:con-a}  is less restrictive than \eqref{e:con-wave}. In fact, for $d\in \{1,2\}$, \eqref{e:con-a} is automatically fulfilled as soon as $a_0\in[0,1)$ and the covariance function $\gamma$ is nonnegative (which implies  $a\in(0,d]$, see \cite[Remark 1.4]{CDST}). Nevertheless, Theorem \ref{thm:wave-eq} gives the first condition on  $a_0$ and $\mu(d\xi)$ in order to get a unique pathwise solution to \eqref{e:wave}. In addition, as the reader will see, our main result also unifies the treatment for $d=1$ and $d=2$.

Let us mention some recent progress, obtained by methods which are completely different from ours, towards the definition of the stochastic wave equation in a pathwise sense.

\noindent
(i) As far as we know, the first pathwise developments for a stochastic wave equation can be found in \cite{qt07}, dealing with the one-dimensional case. Based on the specific expression of the wave kernel for $d=1$, the strategy therein relies on a natural preliminary rotation of the model, that turns \eqref{e:wave} into a more tractable equation in the plane $\R^2$. It is then established that when one injects a noise of the form \eqref{e:cov-sko}-\eqref{e:gamma-scaling} within the new equation (thus corresponding to a \enquote{rotated} noise for the original equation), the interpretation and wellposedness can be obtained for all $a_0,a\in (0,1)$. Unfortunately, this preliminary transform of the model - which will not occur in our direct approach - essentially rules out the possibility to compare the interpretation in \cite{qt07} with ours, and accordingly to compare the conditions in \cite{qt07} with those in \eqref{e:con-wave}.

\noindent (ii) The recent publication \cite{balan22} focuses on the wave equation with a noise $\dot W$ independent of time. The condition obtained in \cite{balan22} is
\begin{equation}\label{e:con-balan}
\int_{\R^d} \frac1{1+|\xi|} \, \mu(d\xi) <\infty,
\end{equation}
which is equivalent to  $a<1$ when $\gamma$ has the scaling property \eqref{e:gamma-scaling}.  A time independent noise being morally equivalent to $a_0=0$, in this situation our condition \eqref{e:con-wave} can be read 
\[a<1 \text{ for } d=1, \quad\text{ and }\quad a<\frac12 \text{ for } d=2. 
\]
Our condition on $a$ is thus slightly suboptimal for $d=2$ if one compares it to~\eqref{e:con-balan}, although our method is currently the only one accommodating for a time-dependent noise.
Let us also mention that \cite{balan22} relies on chaos expansions for Stratonovich integrals. This strategy is interesting in its own right and very different from the pathwise considerations in the current paper. 

\noindent (iii) Some attention has been paid recently to models of wave equation with additive noise and polynomial nonlinearities. A prototype for such an equation can be written as 
\begin{equation}\label{e:wave-additive}
\frac{\partial^2 u}{\partial t^2} (t,x) = \Delta u- \rho u^2 + \dot W,
\end{equation}
with an additive fractional noise $\dot W$  and $\rho>0$. In this context, renormalization procedures are implemented in \cite{deya19,gko18}. This yields existence results for \eqref{e:wave-additive} in case of a rough noise $\dot W$. 

\smallskip

As one can see, the study of wave equations in a rough setting is still wide open. Our contribution aims at a better understanding of the Young regime within this landmark.

As mentioned above, instead of  the Skorohod setting advocated in \cite{CDST}, the stochastic differential in \eqref{e:wave} is interpreted in the {\it Stratonovich} sense.   This forces to a totally different approach which is based on pathwise type considerations.  We briefly elaborate our strategy below. 

We will solve the wave equation under its so-called mild form, which will be properly introduced in \eqref{equa-young-wellposed}. The main technical issue will thus be to understand the meaning of integrals like 
\begin{equation}\label{e:int}
\int_0^t \cg_{t-r}(u_r d\dot W_r), 
\end{equation}
where $u$ is a candidate solution, $\cg_{t}$ stands for the wave operator at time $t$, and where $\dot W$ is our noisy input. Since we wish to interpret \eqref{e:int} as a Young integral, we will see it as the limit of a sequence of paths $\{\cj_t^{(n)}; t\in[0,T]\}$ defined by 
\begin{equation}\label{e:j-n}
\cj_t^{(n)}:=\sum_{k=0}^{m-1} \cg_{t-t^n_k}(u_{t^n_k}\, \der \dot{W}_{t^n_k t^n_{k+1}}), 
\quad\mbox{ for } t\in (t^n_{m-1},t^n_m]\,,
\end{equation}
where $t_m^n = mT/2^n$ for $0\le m \le 2^n$ and $\delta \dot W_{t_k^nt_{k+1}^n}=\dot W_{t_{k+1}^n}-\dot W_{t_k^n}$. Our most important endeavour is thus to quantify the smoothing effect of the operator $\cg_{t-t_k^n}$ in \eqref{e:j-n}. This is achieved thanks to some Strichartz type estimates for the operator $\cg_s$ in some weighted Besov spaces, which are new and should be considered as one of the main contributions of the current paper. 

Our findings concerning Strichartz type estimates are summarized in Proposition \ref{prop:regu-effect} below. Roughly speaking, if $\cb^\alpha$ designates a weighted Besov space  with regularity $\alpha\le 0$, we claim that
\begin{equation}\label{e:strich-est}
\begin{aligned}
\|\{\cg_t-\cg_s\} f\|_{\cb^{\alpha+\ka}} &\lesssim |t-s|^{1-\ka} \|f\|_{\cb^\alpha}, \quad\text{for all} \quad \ka\in [0,1],\ \text{ if } d=1, \\
\text{ and } \|\{\cg_t-\cg_s\} f\|_{\cb^{\alpha+\ka}} &\lesssim |t-s|^{1/2-\ka} \|f\|_{\cb^\alpha}, \quad\text{for all} \quad \ka \in [0,1/2], \ \text{ if } d=2. \\
\end{aligned}
\end{equation}
These bounds thus differ from the classical Strichartz estimates (see e.g. \cite{gv95}) in two aspects. First, they offer a H{\"o}lder-type control (in time) on the action of $\cg$, a key ingredient toward a successful Young integration procedure. Secondly, they involve Besov spaces with weights, which is crucial in order to deal with a noise like $\dot{W}$ defined on the whole space.

Our considerations are based on a weighted version of Littlewood-Paley analysis introduced by Rychkov  in~\cite{rychkov}. This approach is convenient when one wishes to handle kernels with polynomial decay in Fourier modes such as $\cg_t$. This explains  in particular  why we have decided to stick to \cite{rychkov} instead of using the more recent method advocated in \cite{mourrat-weber}.

We close this section by highlighting some possible generalization of our work:

\noindent (i) In the current article, we have focused on a noisy term of the form $u\dot W$ in \eqref{e:wave} for sake of conciseness. We firmly believe that a noise term $\sigma(u) \dot W$ with a smooth enough $\sigma$ could also be covered by our approach, at the price of longer technical considerations.

\noindent (ii)  The treatment of rougher noises is expected to rely on higher-order expansions of the model, together with renormalization procedures. These challenging developments may require the adaptation of some ideas from regularity structures theory or paracontrolled calculus in the wave setting. 

\noindent (iii) The Strichartz type estimate alluded to above are based on convolution estimates for functions. In order to reach $d=3$ or above, those estimates should be extended to measures or distributions. At the moment, this is still an open problem for us.

The rest of this paper is organized as follows.  In Section \ref{sec:besov}, we recall some preliminaries on weighted Besov spaces. The smoothing effect of the wave kernel in weighted Besov spaces is obtained in Section \ref{sec:smoothing-effect}.  It is then applied  to study a Young-type equation associated with \eqref{e:wave} in Section \ref{sec:young-eq}, where we define Young wave integration and prove the wellposedness of the Young-type equations.  In Section \ref{sec:wave}, we apply the result on the Young-type equation obtained in Section \ref{sec:young-eq} to our wave equation \eqref{e:wave} and obtain the existence and uniqueness of the solution under condition \eqref{e:con-wave}.  Finally, in Appendix~\ref{sec:appendix}, we provide  details for some results used in the preceding sections.

 \begin{notation}\label{notations}
 For $r=(r_1, \dots, r_d)$ with $r_i\in \mathbb N$ and $x=(x_1, \dots, x_d)$ with $x_i\in \R^d$, we set  \[|r|=\sum_{i=1}^d r_i, \quad r!=r_1!\dots r_d!, \quad x^r =x_1^{r_1}\dots x_d^{r_d}, \quad D^r f(x) = \left(\frac{\partial}{\partial x_1}\right)^{r_1}\dots \left(\frac{\partial}{\partial x_d}\right)^{r_d}.\]
 For a generic continuous function $\vp$ defined on $\R^d$, we set $\vp_j(x) = 2^{dj} \vp(2^j x).$ We use $\bfb_R(x)$ to denote the open ball in $\R^d$ centered at $x$ with radius $R>0$, i.e., $\bfb_R(x)=\{y\in\R^d; |y-x|<R\}$, and in particular we set $\bfb_R:=\bfb_R(0).$  Finally, we denote by $\cac_c^\infty$ the set of smooth and compactly-supported functions on $\R^d$. 
\end{notation}

\section{Weighted Besov spaces}\label{sec:besov}

Equation \eqref{e:wave} will be solved in some properly weighted Besov spaces.  In fact for a sharp analysis of the wave kernel properties (in Section \ref{sec:smoothing-effect}), we have found convenient to use the general setting introduced by Rychkov in \cite{rychkov}.  This setting has the advantage to cover a large class of weights. The current section is mostly dedicated to recall the main elements of the latter article.

\subsection{Weighted Besov spaces} 
In this section we construct and give some basic properties of the weighted Besov spaces used for the wave equation.  

Let us first recall the definition of the class $A_p^{\textnormal{loc}}$ of weights which was introduced in \cite{rychkov}.
\begin{definition}
We call weight any locally integrable and strictly positive function $w:\R^d\to \R_+$.  Next for  every $1<p<\infty$, we define $A^{\textnormal{loc}}_p$ as the set of weights $w$ for which
\begin{equation}\label{e:weight}
\|w\|_{p,\textnormal{loc}}:=\sup_{|I|\leq 1} \frac{1}{|I|^p} \bigg(\int_I w(x)\, dx\bigg)\bigg(\int_I w(x)^{-\frac{p'}{p}} \, dx\bigg)^{\frac{p}{p'}} \ < \ \infty \, ,
\end{equation}
where $p'$ refers to the H{\"o}lder conjugate of $p$, $I=[a_1, b_1]\times \dots \times[a_d, b_d]$ with $a_i<b_i$ for $i=1,\dots, d$ and $|I|=\prod_{i=1}^d|a_i-b_i|$.
\end{definition}
 As pointed out in \cite{rychkov}, this general class $A_p^{\textnormal{loc}}$ of weights consists of the classical $A_p$ weights (see e.g. \cite{stein93}) and the locally regular weights that grow/decay at most exponentially at infinity. In particular, for $1<p<\infty$, $A_p^{\textnormal{loc}}$ includes the following exponential and polynomial weights that we will use in this paper: 
\begin{equation}\label{e:weights}
w_\mu(x) := e^{-\mu|x|} \quad\text{ for } \mu\in \R , \quad\text{ and }\quad P(x) := \left(1+|x|^{d+1}\right)^{-1}\,. 
\end{equation}

For a generic weight $w$, we  introduce the notion of weighted $L^p$ space which will be at the heart of our analysis.

\begin{definition}\label{def:L-mu}
 Let $w$ be a weight function in the space $A_p^{\textnormal{loc}}$ for $p>1$. For a function $f:\R^d\to\R$, we set 
\[\|f\|_{L^p_w}:=\left(\int_{\R^d} |f(x)|^p \, w(x)\, dx\right)^{1/p}\, .
\] 
Whenever  $w_\mu(x)=e^{-\mu|x|}$ is an  exponential weight, we write $L_{w_\mu}^p:=~L_\mu^p.$

\end{definition}

We also label a notation for functions with centered moments in the definition below. 

\begin{definition}\label{def:D-L}
For every integer $L\geq 0$, we denote by $\cd_L$ the set of functions $\vp\in \cac_c^\infty$ such that for each multi-index $r$ with  $0\leq |r|\leq L$, one has
\begin{equation}\label{e:D-L}
\int_{\R^d} \, x^r \vp(x)\, dx=0 \, .
\end{equation}
\end{definition}

The compatibility of weights with convolution products is a crucial and basic feature for a proper Besov analysis. We state and prove a lemma in this sense.

\begin{lemma}\label{lem:conv-comp}
Fix two constants $  \mu_\ast  ,K>0$. Then for all $0\leq \mu \leq   \mu_\ast  $ and $\vp\in \cac_c^\infty$ such that $\text{Supp}(\vp)\subset \mathbf B_K$, one has

\begin{equation}\label{conv-comp-gene}
\big\| \vp \ast f\big\|_{L^p_\mu}\lesssim \big\| |\vp|\ast \big[ w_\mu^{1/p} |f| \big]\big\|_{L^p}\, ,
\end{equation}
for some proportional constant that only depends on $p$, $K$ and $  \mu_\ast  $. In particular, the following compatibility relation holds true: 
\begin{equation}\label{e:conv-Young}
\big\| \vp \ast f\big\|_{L^p_\mu}%\lesssim \big\| |\vp|\ast \big[ w_\mu(.)^{\frac{1}{p}} |f| \big]\big\|_{L^p}
\lesssim \|\vp\|_{L^1}\|  f\|_{L^p_\mu}.
\end{equation}
\end{lemma}

\begin{proof}
We start by writing the $L_\mu^p$-norm of $\vp\ast f$ according to Definition \ref{def:L-mu}, 
\begin{align*}
\big\| \vp\ast f\big\|_{L^{p}_\mu}^{p}&=\int_{\R^d} dx\, w_\mu(x)\bigg|\int_{\R^d} dy\, \vp(x-y) f(y)\bigg|^{p}. 
\end{align*}

We now insert the weight $w_\mu(y)$ in order to get 
\begin{align}
\big\| \vp\ast f\big\|_{L^{p}_\mu}^{p}&=\int_{\R^d} dx\, \bigg|\int_{\R^d} dy\, \bigg[\bigg(\frac{w_\mu(x)}{w_\mu(y)}\bigg)^{\frac{1}{p}}\vp(x-y)\bigg] \Big[ w_{\mu}(y)^{\frac{1}{p}}f(y)\Big]\bigg|^{p}\notag\\
&\leq \int_{\R^d} dx\, \bigg|\int_{\R^d} dy\, \bigg|\bigg(\frac{w_\mu(x)}{w_\mu(y)}\bigg)^{\frac{1}{p}}\vp(x-y)\bigg| \Big| w_{\mu}(y)^{\frac{1}{p}}f(y)\Big|\bigg|^{p}\, .\label{e:vp*f}
\end{align}
Recall that $\text{Supp}(\vp)$ is a subset of $\mathbf B_K$. Hence for all $x,y\in\R^d$, we have 
\[\bigg|\bigg(\frac{w_\mu(x)}{w_\mu(y)}\bigg)^{\frac{1}{p}}\vp(x-y)\bigg| \le e^{\mu K/p} |\vp(x-y)| .\]
Plugging this inequality into \eqref{e:vp*f} we immediately get \eqref{conv-comp-gene}, which ends our proof. 
\end{proof}

In order to implement the partition of unity necessary for a proper definition of Besov spaces, we state a technical lemma about convolutions of rescaled functions. 

\begin{lemma}\label{lem:vp-j-ast-psi-l}
Fix an integer $L\geq 0$. Let $\vp,\psi\in \cac_c^\infty$ and suppose $\psi\in \cd_L$, where we recall that $\cd_L$ is introduced in Definition \ref{def:D-L}. Next for every $j\geq 1$, recall from Notation \ref{notations} that we have set $\vp_j(x):=2^{dj} \vp(2^j x)$, $\psi_j(x):=2^{dj} \psi(2^j x)$. Then for all $0\leq j\leq \ell$ and $1\leq r\le \infty$, it holds that
\begin{equation}\label{notations-psi-l}
\big\|\vp_j \ast \psi_\ell\big\|_{L^r} \lesssim 2^{dj\left(1-\frac{1}{r}\right)} 2^{-L(\ell-j)}\, .
\end{equation}
\end{lemma}

\begin{proof}
   In order to ease our notation, let us agree on the following convention:

\noindent
\textit{For the sake of clarity,  we  assume from now on that the supports of functions in $\cac_c^\infty$ are all included in the unit ball $\mathbf B_1$.}

\noindent We first prove inequality \eqref{notations-psi-l} for the case $1\le r< \infty$. Writing the definition of the convolution product and performing the change of variables $2^{\ell}y\to y$ and  $2^{j}x\to x$, one gets
\begin{align*}
\big\|\vp_j \ast \psi_\ell\big\|_{L^r}^r &=\int dx\, \bigg| \int dy\, 2^{dj}\vp(2^j(x-y))2^{d\ell}\psi(2^\ell y)\bigg|^r\\
&=2^{dj (r-1)}\int_{\mathbf B_2}dx \, \bigg| \int_{\mathbf B_1} dy \, \psi(y) \vp(x-2^{-(\ell-j)}y)\bigg|^r. 
\end{align*}
Next invoke the fact that $\psi\in \mathcal D_L$ in order to write 
\begin{align}\label{e:phi*psi}
&\big\|\vp_j \ast \psi_\ell\big\|_{L^r}^r  \notag\\
=&2^{dj (r-1)}\int_{\mathbf B_2}dx \, \bigg| \int_{\mathbf B_1} dy \, \psi(y) \bigg[\vp(x-2^{-(\ell-j)}y)-\sum_{0\le |k| \leq L-1} \frac{ D^k\vp(x)}{k!} (-2^{-(\ell-j)}y)^k \bigg]\bigg|^r\, .
\end{align}
Furthermore, thanks to classical considerations on remainders for Taylor expansions,  for $x,y\in \mathbf B_2$, we have
\begin{equation}\label{e:est1}
\bigg|\vp(x-2^{-(\ell-j)}y)-\sum_{0\le |k| \leq L-1} \frac{D^k\vp(x)}{k!} (-2^{-(\ell-j)}y)^k \bigg|\lesssim \, 2^{-L(\ell-j)}.
\end{equation}
Plugging \eqref{e:est1} into \eqref{e:phi*psi}, this proves \eqref{notations-psi-l} for $1\le r<\infty$.

For the case $r=\infty$ in \eqref{notations-psi-l},  we have similarly
 \begin{eqnarray*}
\vp_j \ast \psi_\ell(   2^{-j}x  )&=&  \int dy\, 2^{dj}\vp(2^j(   2^{-j}x  -y))2^{d\ell}\psi(2^\ell y)
=2^{dj} \int_{\mathbf B_1} dy \, \psi(y) \vp(x-2^{-(\ell-j)}y)\\
&=&2^{dj}\int_{\mathbf B_1} dy \, \psi(y) \bigg[\vp(x-2^{-(\ell-j)}y)-\sum_{0\le |k| \leq L-1} \frac{ D^k\vp(x)}{k!} (-2^{-(\ell-j)}y)^k \bigg]\, ,
\end{eqnarray*}
and \eqref{notations-psi-l} follows from \eqref{e:est1}. This finishes the proof.  
\end{proof}

Eventually we will define a bump type function adapted to the construction of our Besov type spaces. 

\begin{notation}\label{not:vp-0-vp}
Pick an arbitrary function $\vp_0\in \cac_c^\infty$. Then we denote by $\vp$ the function defined for $x\in \R^d$ by
$$\vp(x):=\vp_0(x)-2^{-d}\vp_0\Big(\frac{x}{2}\Big)\, .$$
Also recall that, for every continuous function $g:\R^d\to\R$ and $j\geq 1$, we set $g_j(x):=2^{dj} g(2^j x)$ for all $x\in \R$.
\end{notation}

We are now ready to give a precise definition of the Besov spaces considered in this article, which is borrowed from \cite[Definition 2.4]{rychkov}.

\begin{definition}\label{def:Besov} Consider a function $\vp_0\in \cac_c^\infty$  and the corresponding $\vp$ as in Notation~\ref{not:vp-0-vp}. For any fixed length $L$, we assume that $ \int_{\R^d} \vp_0(x)\, dx\neq 0$ and $\vp\in\cd_L$. Otherwise stated we have,
\begin{equation}\label{condition-l-1}\tag{$\ca_\vp$}
\int_{\R^d} \vp_0(x)\, dx\neq 0, \quad \text{and} \quad \int_{\R^d} x^r\vp(x)\, dx= 0\, , \text{ for } |r|\le L\,.
\end{equation}  
Then for all $-\infty<s\leq 1$, all $1<p<\infty$ and $1<q\leq \infty$, we define the space $\cb^{s,\mu}_{p,q}(\vp_0)$ as the completion of $\cac_c^\infty$ with respect to the norm
\begin{equation}\label{e:Besov-norm}
\|f\|_{\cb^{s,\mu}_{p,q}(\vp_0)}:=\bigg(\sum_{j\geq 0} 2^{jsq}\big\|\vp_{j}\ast f\big\|_{L^p_\mu}^q\bigg)^{\frac{1}{q}} \, .
\end{equation}
\end{definition}

At first sight it seems that the spaces $\cb^{s,\mu}_{p,q}(\vp_0)$ depend on the specific function $\vp_0$ we have chosen. The following Proposition, borrowed from \cite[Corollary 2.7]{rychkov}, shows that the Besov spaces do not exhibit this kind of dependence. 

\begin{proposition}\label{prop:equivalence-norms}
Let $\vp_0,\tilde{\vp}_0\in \cac_c^\infty$ be two functions satisfying Assumption \eqref{condition-l-1}, and fix $  \mu_\ast  >0$. Then there exist two constants $c_{\vp_0,\tilde{\vp}_0,  \mu_\ast  },C_{\vp_0,\tilde{\vp}_0,  \mu_\ast  }>0$ such that for all $-\infty <s\leq 1$, all  $0\leq \mu\leq   \mu_\ast  $, and $1<p,q<\infty$ we have 
\begin{equation}
c_{\vp_0,\tilde{\vp}_0,  \mu_\ast  }\|f\|_{\cb^{s,\mu}_{p,q}(\tilde{\vp}_0)}\leq \|f\|_{\cb^{s,\mu}_{p,q}(\vp_0)}\leq C_{\vp_0,\tilde{\vp}_0,  \mu_\ast  }\|f\|_{\cb^{s,\mu}_{p,q}(\tilde{\vp}_0)}\,.
\end{equation}
\end{proposition}

%\color{black}

Since Proposition \ref{prop:equivalence-norms} states that the Besov spaces do not depend on the particular choice of $\varphi_0$, one might just fix a particular $\varphi_0$ and consider the corresponding spaces $\cb$. This convention is labeled below. 

\begin{notation}\label{notation:pmb-vp}
From now on, we fix a function $\pmb{\vp_0}\in \cac_c^\infty$ satisfying \eqref{condition-l-1} and consider the scale of spaces $\cb^{s,\mu}_{p,q}:=\cb^{s,\mu}_{p,q}(\pmb{\vp_0})$, for $1<p<\infty$, $1<q\leq \infty$, $-\infty <s\leq 1$ and $\mu\in \R$.
\end{notation}

\subsection{A few properties of the weighted Besov spaces}

We now collect a few basic facts about Besov spaces, whose role for our existence and uniqueness result will be crucial.

First most of our weights $w$ in the sequel are exponential. However, we also deal occasionally with polynomial weights. The lemma below gives a link between those classes.

\begin{lemma}\label{lem:change-weight}
Let $  \mu_\ast  >0$ be a fixed constant. Consider the polynomial weight $P(x)$ and the family $\{w_\mu, 0< \mu\le   \mu_\ast  \}$ of exponential weights given in \eqref{e:weights}. Then, for $1< p<\infty, 1<q\le \infty$ and $s\in \R$, we have 
$$\big\|f\big\|_{\cb^{s,w_\mu}_{p,q}}\lesssim \mu^{-\frac{d+1}{p}} \big\|f\big\|_{\cb^{s,P}_{p,q}}\, ,$$ 
where the multiplicative constant on the right-hand side depends on $  \mu_\ast  $ only.
\end{lemma}

\begin{proof} Invoking the fact that $0<\mu\le   \mu_\ast  $, it is readily checked that
\[e^{-\mu|x|}\lesssim P(x) \bigg\{ \1_{\{|x|<1\}}+\frac{1}{\mu^{d+1}}\1_{\{|x|\geq 1\}}\bigg\}\lesssim \frac{P(x)}{\mu^{d+1}} \big\{  \mu_\ast^{d+1}  +1\big\}\, .\]
Inserting this inequality into the definition \eqref{e:Besov-norm} of Besov norms, our claim is achieved in a straightforward way. 
\end{proof}

Interpolation inequalities are part of the basic toolkit of Besov spaces. In our weighted context, we can state the following result.

\begin{lemma}\label{lem:interpol}
Consider    nine    parameters $s, s_0, s_1\in(-\infty, 1]$, $p, p_0, p_1\in(1, \infty)$ and $q, q_0, q_1\in(1, \infty]$. We suppose that there exists $\theta\in[0,1]$ such that 
$$s=(1-\theta) s_0+\theta s_1\, , 
\qquad \frac{1}{p}=\frac{1-\theta}{p_0}+\frac{\theta}{p_1}\, , 
\qquad \frac{1}{q}=\frac{1-\theta}{q_0}+\frac{\theta}{q_1}\, .$$
Then    for every $\mu\in \R$,    it holds that
\begin{equation}\label{e:interpol}
\big\| f\big\|_{\cb^{s,\mu}_{p,q}} \leq \big\| f\big\|_{\cb^{s_0,\mu}_{p_0,q_0}}^{1-\theta} \big\| f\big\|_{\cb^{s_1,\mu}_{p_1,q_1}}^\theta \, .
\end{equation}
\end{lemma}

\begin{proof}
The proof is achieved by standard H\"older-type inequalities. It is left to the reader for the sake of conciseness. 
\color{black}
\end{proof}

We next move to an embedding property which will be an important step in our analysis. 

\begin{lemma}\label{lem:embed}
We consider    six    parameters $1<p<p'<\infty$, $-\infty<\ka'<\ka\le 1$ and $0\le \mu< \mu'\le    \mu_\ast  $. Assume that the following relation holds true, 
\begin{equation}\label{e:ka}
 \ka=\ka'+d\bigg(\frac{1}{p}-\frac{1}{p'}\bigg)\, , \quad\text{and}\quad
 \mu'=\frac{p'}{p} \mu\, .
\end{equation}
Then for the Besov spaces introduced in Definition \ref{def:Besov}, it holds that
\begin{equation}
\big\|f \big\|_{\cb^{\ka',\mu'}_{p',\infty}} \lesssim \big\|f \big\|_{\cb^{\ka,\mu}_{p,\infty}}  \, , 
\end{equation}
for some proportional constant that only depends on $  \mu_\ast  $ (as far as weight parameters are concerned). 
\end{lemma}

\begin{proof}
Consider an integer $L\geq 1$, the exact value of which will be determined later on (namely we will see that $L=L(\ka, \ka')$ depends on $(\ka, \ka')$ only). Recall that we have defined our Besov spaces based on a function $\varphi_0$ as in Notation \ref{notation:pmb-vp}.  By Theorem \ref{theo:local-reproducing-formula} in the Appendix, we know that there exist $\psi_0\in \cac_c^\infty$ and $\psi\in \cd_L$ such that
\begin{equation}\label{decompo-besov-embed}
\vp_j \ast f=\sum_{\ell\geq 0}  (\vp_j \ast \psi_\ell) \ast (\vp_\ell\ast f)\, .
\end{equation}

Observe that we have assumed the supports of the $\varphi$ and $\psi$ functions    to be    subsets of $\mathbf B_1$,    and so    we have  $\text{Supp}\, (\vp_j \ast \psi_\ell) \subset \mathbf B_2$ for all $j, \ell \ge 0$. Hence  we can use \eqref{conv-comp-gene} to assert that
$$\big\| (\vp_j \ast \psi_\ell) \ast (\vp_\ell\ast f)\big\|_{L^{p'}_{\mu'}} \lesssim \big\| |\vp_j \ast \psi_\ell| \ast |w_{\mu}^{\frac{1}{p}}(\vp_\ell \ast f)| \big\|_{L^{p'}}\, , $$
for some proportional constant that only depends on $  \mu_\ast  $. We can then invoke the classical Young inequality for convolution products to derive  
\begin{equation}\label{e:convolutions-f}
\big\| (\vp_j \ast \psi_\ell) \ast (\vp_\ell\ast f)\big\|_{L^{p'}_{\mu'}}\lesssim \big\| \vp_j \ast \psi_\ell\big\|_{L^r} \big\|w_{\mu}^{\frac{1}{p}}(\vp_\ell \ast f)\big\|_{L^{p}} \lesssim \big\| \vp_j \ast \psi_\ell\big\|_{L^r} \big\|\vp_\ell\ast f\big\|_{L^{p}_{\mu}}
\end{equation}
where $r$ is defined through the relation $\frac{1}{r}+\frac{1}{p}=1+\frac{1}{p'}$. 
Going back to \eqref{decompo-besov-embed}, we deduce that for every $j\geq 0$,
\begin{align*}
\big\|\vp_j \ast f\big\|_{L^{p'}_{\mu'}}&\lesssim \sum_{\ell\geq 0} \big\| \vp_j \ast \psi_\ell\big\|_{L^r} \big\|\vp_\ell\ast f\big\|_{L^{p}_{\mu}}\lesssim \|f\|_{\cb^{\ka,\mu}_{p,\infty}}\sum_{\ell\geq 0} 2^{-\ka\ell}\big\| \vp_j \ast \psi_\ell\big\|_{L^r} \, .
\end{align*}

We now split the sum above into $\ell \le j$ and $\ell >j$.   We also assume $\vp\in \mathcal D_L$ in Definition~\ref{def:Besov}.   Then we apply Lemma \ref{lem:vp-j-ast-psi-l} in two different ways according to those two cases (recall that \eqref{notations-psi-l} is valid for $j\le \ell$). This yields,  
\begin{align*}
\big\|\vp_j \ast f\big\|_{L^{p'}_{\mu'}}&\lesssim \|f\|_{\cb^{\ka,\mu}_{p,\infty}}\bigg[\sum_{0\leq \ell\leq j} 2^{-\ka \ell}\big\| \vp_j \ast \psi_\ell\big\|_{L^r} + \sum_{\ell>j} 2^{-\ka \ell}\big\| \vp_j \ast \psi_\ell\big\|_{L^r}\bigg]\, , \\
&\lesssim \|f\|_{\cb^{\ka,\mu}_{p,\infty}}\bigg[\sum_{0\leq \ell\leq j} 2^{-\ka \ell}2^{d\ell(1-\frac{1}{r})}2^{-L(j-\ell)} + 2^{dj(1-\frac{1}{r})}\sum_{\ell>j} 2^{-\ka \ell}2^{-L(\ell-j)}\bigg]\,.
\end{align*}
Performing the change of variables $j-\ell\to \ell$ and $\ell-j\to \ell$ in the two sums above, we thus get 
\begin{align*}
\big\|\vp_j \ast f\big\|_{L^{p'}_{\mu'}}&\lesssim \|f\|_{\cb^{\ka,\mu}_{p,\infty}}\bigg[2^{-j(\ka -d(1-\frac{1}{r}))}\sum_{0\leq \ell\leq j} 2^{-(L-\ka+d(1-\frac{1}{r}))\ell} + 2^{-j(\ka -d(1-\frac{1}{r}))}\sum_{\ell>0}2^{-(L+\ka)\ell}\bigg]\, .
\end{align*}
Now remember that $1-\frac{1}{r}=\frac{1}{p}-\frac{1}{p'}$, so that $\ka -d(1-\frac{1}{r})=\ka'$ according to \eqref{e:ka}.  We thus have obtained that
\begin{equation}\label{notations-f}
2^{\ka' j}\big\|\vp_j \ast f\big\|_{L^{p'}_{\mu'}}\lesssim \|f\|_{\cb^{\ka,\mu}_{p,\infty}}\bigg[ \sum_{0\leq \ell\leq j} 2^{-(L-\ka')\ell} +\sum_{\ell>0}2^{-(L+\ka)\ell}\bigg]\, .
\end{equation}
Let us now pick $L>\max(-\ka, \ka')$ in \eqref{notations-f}. We have obtained 
\[   2^{\ka' j}    \big\|\vp_j \ast f\big\|_{L^{p'}_{\mu'}} \lesssim \|f\|_{\cb^{\ka,\mu}_{p,\infty}} \, , \]
where the propotional constant depends on $(L, \ka,\ka')$ only.  To conclude, we simply recall that 
\[\|f\|_{\cb^{\ka',\mu'}_{p',\infty}} =\sup_{j\ge 0}2^{\ka' j} \big\|\vp_j \ast f\big\|_{L^{p'}_{\mu'}}\,,\] 
from which our claim is easily deduced.
\end{proof}

In order to handle the wave equation \eqref{e:wave}, we will have to multiply the candidate solution $u$ with the distribution-valued noise $\dot W$. The lemma below ensures a proper definition of this kind of product in the Young sense. Its proof is deferred to the Appendix (Section~\ref{sec:proof-lem}).

\begin{lemma}\label{lem:multipl}
Fix $  \mu_\ast  >0$. Let  $\alpha, \beta\in(0,1]$, $p, p_1, p_2\in(1,\infty)$ and $\mu,\mu_1,\mu_2\in[0,  \mu_\ast  ]$ be such that 
\begin{equation}\label{e:con-parameters}
\alpha<\beta, \qquad \frac{1}{p}=\frac{1}{p_1}+\frac{1}{p_2},\quad \text{and} \quad \frac{\mu}{p}=\frac{\mu_1}{p_1}+\frac{\mu_2}{p_2} \, . 
\end{equation}
Consider $f\in\cb_{p_1,\infty}^{-\alpha, \mu_1}$ and $g\in \cb_{p_2,\infty}^{\beta, \mu_2}$. Then it holds that
\begin{equation}\label{e:product}
\big\|   f \cdot g  \big\|_{\cb^{-\al,\mu}_{p,\infty}} \lesssim \big\|f \big\|_{\cb^{-\al,\mu_1}_{p_1,\infty}} \big\|g\big\|_{\cb^{\beta,\mu_2}_{p_2,\infty}}  \, ,
\end{equation}
for some proportional constant that only depends on $  \mu_\ast  $.
\end{lemma}

We close this section by rephrasing relation \eqref{e:product} in a way which turns out to be more convenient for our later computations.  

\begin{corollary}\label{cor:useful}
Fix $  \mu_\ast  >0$.  Consider some parameters $\mu,\nu\in[0,   \mu_\ast  ], \alpha, \ka\in(0,1]$ and $p>2$ such that 
\[
   0<\mu-\nu\leq \frac{2\mu_\ast}{p}\,,\qquad \nu\leq \mu_\ast\Big(1-\frac{2}{p}\Big)\, ,   \qquad
\ka>\alpha+\frac dp\,. 
\]
 Then    setting $\varepsilon:=\frac{p}{2}(\mu-\nu)\in (0, \mu_\ast]$  , we have
\begin{equation}\label{e:product-rule}
\big\| f \cdot g\big\|_{\cb^{-\al,\mu}_{2,\infty}}\lesssim \big\|f\big\|_{\cb^{\ka,\nu}_{2,\infty}} \big\|g\big\|_{\cb^{-\al,\varepsilon}_{p,\infty}}\, ,
\end{equation}
for some proportional constant that only depends on $  \mu_\ast  $.
\end{corollary}

\begin{proof}
Set $\bar{\ka}:=\ka-\frac{d}{p}$, and consider $r>2$ such that $\frac1r+\frac1p=\frac12$. We also set  $\bar{\nu}:=\frac{r}{2}\nu\le   \mu_\ast  $.   Since $\bar{\ka}>\al$, we can first apply Lemma \ref{lem:multipl} to assert that
$$
\big\| f \cdot g\big\|_{\cb^{-\al,\mu}_{2,\infty}}\lesssim \big\|f\big\|_{\cb^{\bar{\ka},\bar{\nu}}_{r,\infty}} \big\|g\big\|_{\cb^{-\al,\varepsilon}_{p,\infty}}\, .
$$
It now only remains to apply Lemma \ref{lem:embed}, which yields $\|f\|_{\cb^{\bar{\ka},\bar{\nu}}_{r,\infty}} \lesssim \|f\|_{\cb^{\ka,\nu}_{2,\infty}}$.
\end{proof}

\section{Smoothing effect of the wave kernel}
\label{sec:smoothing-effect}

%\begin{lemma}\label{lem:vp-0}
%For all $t\in [0,1]$, $1\leq p<\infty$ and $-\infty <\mu\leq \mu_0$, it holds that
%$$\| \cg_t f\|_{L^p_\mu} \lesssim \|f\|_{L^p_\mu}\, ,$$
%where the proportional constant only depends on $\mu_0$ and $p$. As a consequence, for every $\al\in \R$, 
%\begin{equation}\label{no-regu}
%\|\cg_t f\|_{\cb^{\al,\mu}_{p,q}} \lesssim \|f\|_{\cb^{\al,\mu}_{p,q}}\, ,
%\end{equation}
%where the proportional constant only depends on $\mu_0$ and $p$. 
%\end{lemma}

%\begin{proof}
%One has
%\begin{align*}
%\| \cg_t f\|_{L^p_\mu}^p &=\int_{\R} dx \, w_\mu(x) \bigg| \int_{\R} dy \, \cg_t(x-y) f(y)\bigg|^p=\int_{\R} dx \, \bigg| \int_{\R} dy \, \bigg(\frac{w_\mu(x)}{w_\mu(y)}\bigg)^{\frac{1}{p}}\cg_t(x-y) \big[ w(y)^{\frac{1}{p}}f(y)\big]\bigg|^p\, .
%\end{align*}
%Since $\text{Supp}\, \cg_t \subset [-1,1]$, we deduce that
%\begin{align*}
%\| \cg_t f\|_{L^p_\mu}^p &\lesssim \int_{\R} dx \, \bigg| \int_{\R} dy \,|\cg_t(x-y)| \big| w_\mu(y)^{\frac{1}{p}}f(y)\big|\bigg|^p\, ,
%\end{align*}
%where the proportional constant only depends on $\mu_0$ and $p$.

%\smallskip

%We can now apply the classical Young inequality for convolution and assert that
%$$\| \cg_t f\|_{L^p_\mu}^p \lesssim \| \cg_t\|_{L^1}^p \|f\|_{L^p_\mu}^p\lesssim  \|f\|_{L^p_\mu}^p \, .$$

%\

%As for \eqref{no-regu}, it now suffices to observe that for every $\ell\geq 0$,
%$$2^{\ell\al}\| \pmb{\vp}_\ell \ast\cg_t f\|_{L^p_\mu}=2^{\ell\al}\|\cg_t( \pmb{\vp}_\ell\ast f)\|_{L^p_\mu}\lesssim 2^{\ell\al}\|\pmb{\vp}_\ell\ast f\|_{L^p_\mu}\, .$$

%\end{proof}

The wave kernel plays a fundamental role in the mild formulation of \eqref{e:wave}. In this section, we investigate Strichartz type inequalities for this kernel in the weighted Besov spaces defined in Section \ref{sec:besov}. 
 
\subsection{Definition of the wave kernel}\label{sec:wave-kernel}

In the sequel, we will write $\cg_t$ for the wave operator on $\R^d$, generated by $\partial^2_{tt}-\Delta$. The corresponding kernel will be denoted by $G_t$. For $d=1, 2$, the kernel $G_t$ has the following explicit expression:

\begin{equation}\label{e:G}
G_t(x)=\left\{\begin{array}{lll}
\displaystyle \frac12\1_{[|x|<t]}& \text{ if } d=1,\vspace{0.2cm}\\
\displaystyle\frac1{2\pi} \frac1{\sqrt{t^2-|x|^2}}\1_{[|x|<t]} & \text{ if } d=2.
\end{array} \right.
\end{equation}

It is often easier to express $G_t$ in Fourier modes. Namely for a function $g(t,x)$ we set $\cf g(t,\xi)$ or $\hat g(t, \xi)$ for the spatial Fourier transform, defined by 
\[\hat g(t,\xi) =\int_{\R^d} e^{-\imath x\cdot \xi} g(t,x) dt\,. \]
Then the expression for the Fourier transform of $G_t$ is 
\begin{equation}\label{e:Fourier-G}
\hat G_t(\xi) =\frac{\sin(t|\xi|)}{|\xi|}\,, \xi\in\R^d\,.
\end{equation}

\subsection{Strichartz type estimates}

The aim of this section is to quantify the smoothing effects of the operator introduced in Section \ref{sec:wave-kernel}. We start by  two propositions  giving some information about the behavior in time. Note that we are restricted to a spatial dimension $d=1,2$ for this section.  We  first state  a regularity estimate in dimension 1. 

\begin{proposition}\label{prop:lipschitz}
Assume $d=1$. Let $   \mu_\ast  >0$ be a given constant, and recall that $\cg_t$ designates the wave operator.  Consider $\mu\in[0,    \mu_\ast  ]$ and $0\le s<t\le 1$. Then, for all $f\in L_{\mu}^p$, we have
\begin{equation}\label{e:G-Lip}
\| \{\cg_t-\cg_s\} f\|_{L^p_\mu} \lesssim |t-s|\|f\|_{L^p_\mu}\, ,
\end{equation}
where the proportional constant only depends on $   \mu_\ast  $ and $p$. As a consequence, for all $-\infty <\al \leq 1$ and $1<p, q\leq \infty$, 
\begin{equation}\label{lipschitz}
\|\{\cg_t-\cg_s\} f\|_{\cb^{\al,\mu}_{p,q}} \lesssim |t-s| \|f\|_{\cb^{\al,\mu}_{p,q}}\, ,
\end{equation}
where the proportional constant only depends on $   \mu_\ast  $, $p$ and $q$. 
\end{proposition}

\begin{proof}
As in the proof Lemma \ref{lem:conv-comp}, we write Definition \ref{def:L-mu} for the $L^p$ norm and recall that $\cg_t$ admits the kernel $G_t$ given by  \eqref{e:G}. We get 
\[\| \{\cg_t-\cg_s\} f\|_{L^p_\mu}^p =\int_{\R} dx \, w_\mu(x) \bigg| \int_{\R} dy \, \{G_t-G_s\}(x-y) f(y)\bigg|^p .
\]
Still like in Lemma \ref{lem:conv-comp}, we insert the weight $w_\mu(y)$ and invoke the fact that $\text{Supp}\,(G_t-
G_s)\subset [-2, 2]$. 
This yields
\begin{align*}
\| \{\cg_t-\cg_s\} f\|_{L^p_\mu}^p 
&=\int_{\R} dx \, \bigg| \int_{\R} dy \, \bigg(\frac{w_\mu(x)}{w_\mu(y)}\bigg)^{\frac{1}{p}}\{G_t-G_s\}(x-y) \big[ w_\mu(y)^{\frac{1}{p}}f(y)\big]\bigg|^p\\
&\lesssim \int_{\R} dx \, \bigg| \int_{\R} dy \, \big|\{G_t-G_s\}(x-y)\big| \big| w_\mu(y)^{\frac{1}{p}}f(y)\big|\bigg|^p\, ,
\end{align*}
where the proportional constant only depends on $  \mu_\ast  $ and $p$. 

\smallskip

Taking into account the expression \eqref{e:G}
for $d=1$, it is readily checked that 
\begin{equation}\label{e:G-Lip'}
\|G_t-G_s\|_{L^1(\R)}\le |t-s|\,.
\end{equation}
Therefore, a direct application of Young's inequality (see also relation \eqref{e:conv-Young}) leads to 
\begin{equation}\label{appli-young}
\| \{\cg_t-\cg_s\} f\|_{L^p_\mu}^p \lesssim \|G_t-G_s\|_{L^1}^p \|f\|_{L^p_\mu}^p\lesssim  |t-s|^p\|f\|_{L^p_\mu}^p \, ,
\end{equation}
   which    proves our claim \eqref{e:G-Lip}.

As for \eqref{lipschitz}, it now suffices to observe that for every $\ell\geq 0$,
$$2^{\ell\al}\| \pmb{\vp}_\ell \ast \{\cg_t-\cg_s\}f\|_{L^p_\mu}=2^{\ell\al}\|\{\cg_t-\cg_s\}( \pmb{\vp}_\ell\ast f)\|_{L^p_\mu}\lesssim 2^{\ell\al}|t-s|\|\pmb{\vp}_\ell\ast f\|_{L^p_\mu}\, .$$
Inserting this expression into the definition \eqref{e:Besov-norm} of Besov norm, the desired inequality is achieved. 
\end{proof}

We now turn to  the corresponding property  in the 2-dimensional case. We start with a lemma which summarizes the distinction between the 1-d and the 2-d cases. To get a quick understanding of this, recall from \eqref{e:G-Lip'} that for $d=1$
\[ \|G_t-G_s\|_{L^1(\R)}\le |t-s|\,,\]
while the same result for $d=2$ reads

\begin{lemma} Assume $d=2$ and consider the kernel $G_t$ defined by \eqref{e:G}. Then for all $0\leq s<t\leq 1$,  we have
\begin{equation}\label{e:G-Lip'2}
\big\|G_t-G_s\big\|_{L^1} \lesssim |t-s|^{\frac12}\, .
\end{equation}
\end{lemma}

\begin{proof}
According to relation \eqref{e:G}, for $d=2$ we have
\begin{align*}
\big\|G_t-G_s\big\|_{L^1}&=\frac{1}{2\pi} \int_{\R^2} dx \, \bigg|\frac{1}{\sqrt{t^2-|x|^2}} \1_{\{[x|<t\}}-\frac{1}{\sqrt{s^2-|x|^2}} \1_{\{[x|<s\}}\bigg|\\
&\lesssim \int_{s<|x|<t}\frac{dx}{\sqrt{t^2-|x|^2}}+ \int_{|x|<s}dx \,\bigg|\frac{1}{\sqrt{t^2-|x|^2}}-\frac{1}{\sqrt{s^2-|x|^2}}\bigg|\,. 
\end{align*}
Next we use a polar change of coordinates and set $r:=r^2$ in order to get 
\begin{align*}
\big\|G_t-G_s\big\|_{L^1}
&\lesssim \int_s^t \frac{r\, dr}{\sqrt{t^2-r^2}}+\int_0^s dr \, r\,\bigg|\frac{1}{\sqrt{t^2-r^2}}-\frac{1}{\sqrt{s^2-r^2}}\bigg|\\
&\lesssim \int_{s^2}^{t^2} \frac{dr}{\sqrt{t^2-r}}+\int_0^{s^2} dr \,\bigg|\frac{1}{\sqrt{t^2-r}}-\frac{1}{\sqrt{s^2-r}}\bigg|\, .
\end{align*}
Furthermore, since $s<t$ we have $(s^2-r)^{-1/2}>(t^2-r)^{-1/2}$ for $r\in[0, s^2]$.  This yields 
\begin{align*}
\big\|G_t-G_s\big\|_{L^1}
&\lesssim \int_{s^2}^{t^2} \frac{dr}{\sqrt{t^2-r}}+\int_0^{s^2} dr \,\bigg(\frac{1}{\sqrt{s^2-r}}-\frac{1}{\sqrt{t^2-r}}\bigg)\,.
\end{align*}
We now simply integrate the right-hand side above in order to get (recall that we assume $s,t\in[0,1]$ for this lemma) 
\[ \|G_t-G_s\|_{L^1} \lesssim \Big((t+s)^{\frac12}-(t-s)^{\frac12} \Big) (t-s)^{\frac12}\lesssim (t-s)^{\frac12}\,,\]
which proves our claim \eqref{e:G-Lip'2}. 
\end{proof} 

As a direct consequence, we have the following result parallel to Proposition \ref{prop:lipschitz}. 

\begin{proposition}\label{prop:lipschitz2}
Assume $d=2$. Let $   \mu_\ast  >0$ be a given constant, and recall that $\cg_t$ designates the wave operator.  Consider $\mu\in[0,    \mu_\ast  ]$ and $0\le s<t\le 1$. Then, for all $f\in L_{\mu}^p$, we have
\begin{equation}\label{e:G-Lip2}
\| \{\cg_t-\cg_s\} f\|_{L^p_\mu} \lesssim |t-s|^{   \frac12  }\|f\|_{L^p_\mu}\, ,
\end{equation}
where the proportional constant only depends on $   \mu_\ast  $ and $p$. As a consequence, for all $-\infty <\al \leq 1$ and $1<p, q\leq \infty$, 
\begin{equation}\label{lipschitz2}
\|\{\cg_t-\cg_s\} f\|_{\cb^{\al,\mu}_{p,q}} \lesssim |t-s|^{  \frac12  }\|f\|_{\cb^{\al,\mu}_{p,q}}\, ,
\end{equation}
where the proportional constant only depends on $   \mu_\ast  $, $p$ and $q$. 
\end{proposition}

\begin{proof}
The proof follows exactly that of Proposition \ref{prop:lipschitz}, except for replacing \eqref{e:G-Lip'} by relation~\eqref{e:G-Lip'2}. 
\end{proof}

We now present our main result for this section, quantifying the smoothing effect of $\cg_t$ in our Besov scale. 

\begin{proposition}\label{prop:regu-effect}
In this proposition the dimension parameter takes the values $1$ or $2$.  As in Proposition \ref{prop:lipschitz}, we consider a constant $   \mu_\ast  >0$ and the wave operator $\cg_t$. Let $\alpha, \mu, p,q$ be four parameters with 
\begin{equation}\label{e:4-parameters}
-\infty<\alpha\le 0, 
\qquad
 0\le \mu\le    \mu_\ast  , 
 \qquad
 1<p<\infty, 
 \qquad
 1<q\le \infty.
\end{equation}
Recall that the weighted Besov spaces are introduced in Definition \ref{def:Besov}. Then for all $t\in[0,1]$, it holds that 
\begin{equation}\label{regu-effect-wave}
\big\|\cg_t f\big\|_{\cb_{p,q}^{\al+\rho_d,\mu}} \lesssim \|f\|_{\cb_{p,q}^{\al,\mu}}\, , \text{ with } \rho_d:=\begin{cases}
1 & \text{if} \ d=1\\
  \frac12   & \text{if} \ d=2
\end{cases}\,,
\end{equation}
and where the proportional constant only depends on $p$, $q$ and $   \mu_\ast  $.
\end{proposition}

Before starting the proof of Proposition~\ref{prop:regu-effect}, we label a technical lemma about integrals of the wave kernel.

\begin{lemma}\label{lem:cj-d}
In the setting of Proposition \ref{prop:regu-effect}, and for all $t\in [0,1]$, $j\geq 0$, $y\in \R^d$, define the quantity
\begin{equation}\label{quant-cj}
\ck^{(d)}_{t,j}(y):=2^{j\rho_d}\int_{\bfb_2} dz \, \big|G_t(y-2^{-j}z)-G_t(y)\big|\, .
\end{equation}
Then it holds that
\begin{equation}\label{asser-cj}
\text{Supp}\, \ck^{(d)}_{t,j} \subset \bfb_{t+2} , \quad \text{and} \quad  \sup_{j\geq 0}\sup_{t\in [0,1]} \big\|\ck^{(d)}_{t,j}\big\|_{L^1} \, < \infty.
\end{equation}
\end{lemma}

\begin{proof}
See Appendix.
\end{proof}

\begin{remark}
The presence of a $2^{j\rho_d}$ term in front of the integral in \eqref{quant-cj} is what will guarantee a smoothing effect from $\cb^\alpha$ to a $\cb^{\alpha+\rho_d}$ space. 
\end{remark}

\begin{proof}[Proof of Proposition \ref{prop:regu-effect}]
We begin this proof like in Lemma \ref{lem:embed}. Namely we introduce a parameter $L\ge 0$ whose exact value is to be determined later on. Then we assume that the construction of our weighted Besov spaces hinges on a function $\vp_0$ such that $\vp\in \mathcal D_L$ (where $\vp$ is defined along Notation \ref{not:vp-0-vp}).

Thanks to Theorem \ref{theo:local-reproducing-formula} in the Appendix, we know that there exist
two functions $\psi_0,\psi\in \cac_c^\infty$ such that $\psi\in \mathcal D_L$ and for every $\ell\ge 1$
we have
\begin{equation}\label{decomp-start}
\vp_\ell \ast\cg_t f=\vp_\ell\ast G_t\ast f=\sum_{j\geq 0} (G_t\ast \vp_\ell\ast \psi_j)\ast (\vp_j \ast f) \, ,
\end{equation}
where we recall that $G_t$ is the wave kernel given by \eqref{e:G}. 
In relation \eqref{decomp-start}, recall that we have assumed the supports of $\vp_0$, $\vp$, $\psi$ and $\psi_0$ to be all included in the unit ball $\bfb_1$. However, the arguments could be easily extended to any supporting ball $\bfb_K$.
Let us now split the above sum \eqref{decomp-start} into 
\begin{multline}\label{defi-i-ii}
\sum_{j\geq 0} (G_t\ast \vp_\ell\ast \psi_j)\ast (\vp_j \ast f)
=\sum_{0\leq j< \ell} (G_t\ast \vp_\ell\ast \psi_j)\ast (\vp_j \ast f)
+(G_t\ast \vp_\ell\ast \psi_\ell)\ast (\vp_\ell \ast f)  \\
+\sum_{j>\ell} (G_t\ast \vp_\ell\ast \psi_j)\ast (\vp_j \ast f)
:=I_\ell+I\! I_\ell+I\! I\! I_\ell \, .
\end{multline}
Our inequality \eqref{regu-effect-wave} is now easily seen to be reduced to prove 
\begin{equation}\label{regu-effect-wave'}
\bigg( \sum_{\ell \ge 0} 2^{\ell (\alpha+\rho_d) q} \left(\|I_\ell\|^q_{L_\mu^p}+\|I\!I_\ell\|^q_{L_\mu^p}+\|I\!I\!I_\ell\|^q_{L_\mu^p} \right)\bigg)^{1/q}\lesssim \|f\|_{\cb_{p,q}^{\al,\mu}}\,.
\end{equation}
We will handle the 3 terms on the left-hand side of \eqref{regu-effect-wave'} separately.

\noindent
\textit{Step 1: Bound on $I_\ell$}.  With \eqref{regu-effect-wave} and \eqref{regu-effect-wave'} in mind, our aim is to prove an inequality of the form
\begin{equation}\label{e:I}
\bigg(\sum_{\ell \ge 0} 2^{\ell(\alpha+\rho_d) q}\|I_\ell\|_{L_\mu^p}^q\bigg)^{\frac1q} \lesssim \|f\|_{\cb_{p,q}^{\alpha, \mu}}\,.
\end{equation}
Our first claim is that for $0\le j<\ell$ we can write 
\begin{equation}\label{e:vp-psi}
\vp_\ell\ast \psi_j=(\bar{\psi}\ast \vp_{\ell-j})_j\, , \quad \text{where} \quad \bar{\psi}:=
\begin{cases}
\psi_0 & \text{if} \ j=0\\
\psi &\text{if} \ j\geq 1
\end{cases}\, .
\end{equation}
Indeed, relation \eqref{e:vp-psi} is trivial whenever $j=0$ (note that we have used the convention $g_0=g$ for any given $g\in \cac_c^\infty$ in the right-hand side of \eqref{e:vp-psi}). For $j\ge 1$, the left-hand side of \eqref{e:vp-psi} can be expressed as 
\begin{equation}\label{e:l-vp-psi}
(\vp_\ell\ast \psi_j)(x) = 2^{\ell+j}\int_{\R^d} \vp(2^\ell (x-y)) \psi (2^j y) dy\,,
\end{equation}
while the right-hand side of \eqref{e:vp-psi} has the form
\begin{equation}\label{e:r-vp-psi}
(\vp_{\ell-j}\ast \psi)_j(x) = 2^{\ell} \int_{\R^d} \vp(2^\ell x -2^{\ell-j}y) \psi (y) dy\,.
\end{equation}
Then an elementary change of variable $z=2^{-j}y$ shows that \eqref{e:l-vp-psi} and \eqref{e:r-vp-psi} are equal.

Invoking \eqref{e:vp-psi}, together with the change of variable $2^j z\to z$, it is then readily checked that for $0\le j< \ell$ and $y\in\R$ the term $G_t*\varphi_\ell*\psi_j$ in \eqref{defi-i-ii} can be decomposed as
\begin{equation}\label{fir}
(G_t\ast \vp_\ell\ast \psi_j)(y)=\int_{\R^d} dz \, (\bar{\psi}\ast \vp_{\ell-j})_j(z) G_t(y-z)=\int_{\R^d} dz \, (\bar{\psi}\ast \vp_{\ell-j})(z) G_t(y-2^{-j}z)\, .
\end{equation}
Let us take advantage of some cancellations in the right-hand of \eqref{fir}. Indeed, owing to the fact that $\int_{\R^d} \vp_k(x) dx=0$ for all $k>0$, we also have   
\[G_t(y)\int_{\R^d} dz \, (\bar{\psi}\ast \vp_{k})(z)=0.\]
Subtracting the above quantity to the right-hand side of \eqref{fir}, we get that for all $\ell >j$, 
\begin{align*}
(G_t\ast \vp_\ell\ast \psi_j)(y)&=\int_{\R^d} dz \, (\bar{\psi}\ast \vp_{\ell-j})(z) \big\{G_t(y-2^{-j}z)-G_t(y)\big\} .\end{align*}
In addition,  observing $\text{Supp} \, (\bar{\psi}\ast \vp_{\ell-j}) \subset \bfb_2$, we have that 
%the term $G_t\ast \vp_\ell\ast \psi_j$ in 
\begin{align}\label{g-l-ast'}
(G_t\ast \vp_\ell\ast \psi_j)(y)&=\int_{\bfb_2} dz \, (\bar{\psi}\ast \vp_{\ell-j})(z) \big\{G_t(y-2^{-j}z)-G_t(y)\big\} \, .
\end{align}
With \eqref{g-l-ast'} in hand, the key point now is that due to \eqref{e:D-L} one has easily, for all $0\leq j<\ell$,
\begin{equation}\label{bound-vp-j-ast-psi}
\|\bar{\psi}\ast \vp_{\ell-j}\|_\infty \lesssim 2^{-L(\ell-j)} \, ,
\end{equation}
which follows from Lemma \ref{lem:vp-j-ast-psi-l} directly. Therefore, using the notation introduced in~\eqref{quant-cj}, we have obtained that 
\begin{equation}\label{e:G*vp*psi}
|(G_t*\vp_\ell*\psi_j)(y)| \le 2^{-L(\ell-j)}2^{-j\rho_d}\ck^{(d)}_{t,j}(y)\, .
\end{equation}
Let us insert \eqref{e:G*vp*psi} in $I_\ell$ in order to upper bound this term. That is according to the definition \eqref{defi-i-ii} of $I_\ell$, it is obvious that  
\begin{align}
|I_\ell(x)|&\leq \sum_{0\leq j< \ell} \int_{\R^d} dy \, |(G_t\ast \vp_\ell\ast \psi_j)(y)| |(\vp_j \ast f)(x-y)|\nonumber\\
&\leq \sum_{0\leq j< \ell}2^{-L(\ell-j)}2^{-j\rho_d} \int_{\R} dy \, \ck^{(d)}_{t,j}(y) |(\vp_j \ast f)(x-y)|\nonumber\\
&\leq \sum_{0\leq j< \ell}2^{-L(\ell-j)}2^{-j\rho_d} \big[ \ck^{(d)}_{t,j} \ast |\vp_j \ast f|\big](x).\label{e:I-l'}
\end{align}

We are now ready to analyze the weighted $L^p$ norms of the function $I_{\ell}$. To this aim, we multiply \eqref{e:I-l'} by $2^{\ell(\alpha+\rho_d)}$ and compute the $L^p_\mu$-norm given in Definition \ref{def:L-mu} for $\mu\le   \mu_\ast  $. We get
\begin{align}
2^{\ell(\al+\rho_d)}\|I_\ell\|_{L^p_\mu}&\lesssim \sum_{0\leq j< \ell} 2^{\ell(\al+\rho_d)}2^{-L(\ell-j)}2^{-j\rho_d }2^{-j\al} \Big[ 2^{j\al}\big\|\ck^{(d)}_{t,j} \ast |\vp_j \ast f|\big\|_{L^p_\mu}\Big] \nonumber\\
&\lesssim \sum_{0\leq j< \ell} 2^{(\ell-j)(\al+\rho_d-L)} \Big[ 2^{j\al}\big\|\ck^{(d)}_{t,j} \ast |\vp_j \ast f|\big\|_{L^p_\mu}\Big] .\label{simp-proof-1}
\end{align}
At this point, we can combine the results of Lemma \ref{lem:conv-comp}, Lemma \ref{lem:cj-d}  and Young's inequality to assert that
\begin{align}\label{simp-proof-2}
\big\|\ck^{(d)}_{t,j} \ast |\vp_j \ast f|\big\|_{L^p_\mu}&\lesssim  \big\| |\ck^{(d)}_{t,j}| \ast w_\mu^{1/p} |\vp_j \ast f|\big\|_{L^p}\notag\\ 
&\le \big \| \ck^{(d)}_{t,j} \big\|_{L^1}  \big\| w_\mu^{1/p} |\vp_j \ast f| \big\|_{L^p}\notag\\
&\lesssim \big\|\vp_j \ast f\big\|_{L^p_\mu} \, ,
\end{align}
for some proportional constant that depends neither on $j$ nor on $t$. Injecting this uniform bound into \eqref{simp-proof-1}, we deduce
\begin{equation}\label{e:I-l-Lp'}
2^{\ell(\al+\rho_d)}\|I_\ell\|_{L^p_\mu}\lesssim \sum_{0\leq j< \ell} 2^{(\ell-j)(\al+\rho_d-L)} \Big[ 2^{j\al}\big\|\vp_j \ast f\big\|_{L^p_\mu}\Big]\, .
\end{equation}

We shall recast the right-hand side of \eqref{e:I-l-Lp'} as a discrete convolution. Namely, if we set 
\[
u_j:=2^{j(\al+\rho_d-L)}\1_{\{j\geq 1\}} \, , \quad \text{and}\quad v_j:=2^{j\al}\|\vp_j \ast f\|_{L^p_\mu} \1_{\{j\geq 0\}}\, ,
\]
and if $\ast_{\mathbb{Z}}$ refers to the discrete convolution in $\mathbb Z$, then \eqref{e:I-l-Lp'} can be written as 
\begin{equation}\label{e:I-l-Lp-conv}
2^{\ell(\al+\rho_d)}\|I_\ell\|_{L^p_\mu} \lesssim (u \ast_{\mathbb{Z}} v)_\ell\,.
\end{equation}
The inequality above enables the application of Young's inequality, similarly to what we did in \eqref{e:convolutions-f}. We end up with 
\begin{align}\label{bound-i-l'}
\big\|\ell \mapsto 2^{\ell(\al+\rho_d)}\|I_\ell\|_{L^p_\mu}\big\|_{\ell^q(\mathbb{N})}&\lesssim \big\|j\mapsto 2^{j(\al+\rho_d-L)}\big\|_{\ell^1(\mathbb{N})} \, \big\|j\mapsto 2^{j\al}\|\vp_j \ast f\|_{L^p_\mu}\big\|_{\ell^q(\mathbb{N})}\,. 
\end{align}
Let us finally pick $L >\rho_d+\alpha$, so that the above estimate yields
\begin{equation}\label{e:est-I}
\bigg(\sum_{\ell\ge0} 2^{\ell(\alpha+\rho_d )q}\|I_\ell\|^q_{L_\mu^p}\bigg)^{1/q}\lesssim \|f\|_{\cb^{\alpha,\mu}_{p,q}}\,,
\end{equation}
which is our claim \eqref{e:I} and  the first part of relation \eqref{regu-effect-wave'}.

\noindent
\textit{Step 2: Bound on $I\! I_\ell$}. Recall from \eqref{defi-i-ii} that $I\!I_\ell=(G_t*\vp_\ell*\psi_\ell)*(\vp_\ell* f)$. Our aim for this step is to obtain inequality \eqref{e:est-I} for the sequence $\{I\!I_\ell, \ell \ge0\}$. Whenever $\ell=0$, we appeal again to Lemma \ref{lem:conv-comp} in order to write 
\begin{align*}
\| I\! I_0\|_{L^p_\mu}&\lesssim \|G_t\ast \vp_0\ast \psi_0\|_{L^1} \|\vp_0\ast f\|_{L^p_\mu}\lesssim \|G_t\|_{L^1} \|\vp_0\|_{L^1}\|\psi_0\|_{L^1} \|\vp_0\ast f\|_{L^p_\mu}\,.
\end{align*}
Furthermore, $ \|\vp_0\|_{L^1}, \|\psi_0\|_{L^1} $ above are finite quantities. Since $t\le 1$, the norm $\|G_t\|_{L^1} $ is also bounded, along the same lines as for \eqref{e:G-Lip}. We thus get 
\begin{align}\label{e:II0}
\| I\! I_0\|_{L^p_\mu}\lesssim \|\vp_0\ast f\|_{L^p_\mu} \, ,
\end{align}
where the proportional constant only depends on $p$ and $  \mu_\ast  $.

For $\ell \geq 1$, observe that similarly to \eqref{e:vp-psi} we have $\vp_\ell \ast \psi_\ell=(\vp\ast \psi)_\ell$. We can then use the same arguments as in \eqref{fir}-\eqref{e:G*vp*psi} (replace $\vp_{\ell-j}$ with $\vp$, $\bar{\psi}$ with $\psi$ and $j$ with $\ell$) to derive that for all $y\in \R$,
$$\big|(G_t\ast \vp_\ell\ast \psi_\ell)(y)\big|\lesssim \|\vp\ast \psi\|_\infty 2^{-\ell \rho_d}\ck^{(d)}_{t,\ell}(y)\, .$$
Just as in Step 1 (see \eqref{simp-proof-1}-\eqref{simp-proof-2}), this entails that for $\ell\geq 1$,
$$2^{\ell(\al+\rho_d)}\|I\!I_\ell\|_{L^p_\mu} \lesssim 2^{\ell\al}\|\ck^{(d)}_{t,\ell}\|_{L^1}\|\vp_{\ell}\ast f\|_{L^p_\mu}\lesssim 2^{\ell\al}\|\vp_{\ell}\ast f\|_{L^p_\mu}$$
for some proportional constant that only depends on $p$ and $  \mu_\ast  $. Summing the above inequality over $\ell$ and applying Definition \ref{def:Besov}, we have shown the second part of  relation~\eqref{regu-effect-wave'}, that is
\begin{equation}\label{e:est-II}
\bigg(\sum_{\ell\ge0} 2^{\ell(\alpha+\rho_d)q}\|I\!I_\ell\|^q_{L_\mu^p}\bigg)^{1/q}\lesssim \bigg( \sum_{\ell\ge0} 2^{\ell\alpha q} \| \vp_{\ell}\ast f\|^q_{L_\mu^p} \bigg)^{1/q}\lesssim \|f\|_{\cb_{p,q}^{\alpha, \mu}}\,.
\end{equation}

\noindent
\textit{Step 3: Bound on $I\! I\! I_\ell$}.  We now treat the term $I\!I\!I_\ell$ in \eqref{defi-i-ii}. This will be done along the same lines as for the other steps, and we thus skip some details of computation. First, since $I\!I\!I_\ell$ involves a sum over $j>\ell$, we replace \eqref{e:vp-psi} by 
$$\vp_\ell\ast \psi_j=(\bar{\vp}\ast \psi_{j-\ell})_\ell\, , \ \text{where} \ \bar{\vp}:=
\begin{cases}
\vp_0 & \text{if} \ \ell=0\\
\vp &\text{if} \ \ell\geq 1
\end{cases}\, .$$
Thus, for all $j> \ell$ and $y\in \R$, the equivalent of \eqref{fir} is 
$$
(G_t\ast \vp_\ell\ast \psi_j)(y)=\int_{\R^d} dz \, (\bar{\vp}\ast \psi_{j-\ell})_\ell(z)G_t(y-z)=\int_{\R^d} dz \, (\bar{\vp}\ast \psi_{j-\ell})(z) G_t(y-2^{-\ell}z)\, .
$$

Just as in Step 1, it holds that $\int_{\R} dz \, (\bar{\vp}\ast \psi_{j-\ell})(z)=0$ and $\text{Supp} \, (\bar{\vp}\ast \psi_{j-\ell})\subset \bfb_2$ for all $j>\ell$. Accordingly, relation \eqref{g-l-ast'} for the case $j>\ell$ becomes
\begin{align*}
(G_t\ast \vp_\ell\ast \psi_j)(y)&=\int_{\bfb_2} dz \,  (\bar{\vp}\ast \psi_{j-\ell})(z)  \big\{G_t(y-2^{-\ell}z)-G_t(y)\big\} \, .
%&=\int_{|z|<\inf(2^{\ell+1}t,2)} dz \, (\vp\ast \psi_{j-\ell})(z) \big\{\1_{\{|y-2^{-\ell}z|\leq t\}}-\1_{\{|y|\leq t\}}\big\} \nonumber\\
%&\hspace{1cm}+\int_{\inf(2^{\ell+1}t,2)<|z|<2} dz \, (\vp\ast \psi_{j-\ell})(z)  \big\{\1_{\{|y-2^{-\ell}z|\leq t\}}-\1_{\{|y|\leq t\}}\big\} \, .\label{g-l-ast}
\end{align*}
Next using the notation in \eqref{quant-cj}, we obtain that for all $j> \ell$ and $y\in \R$,
\begin{equation*}%\label{g-t-ast-vp-bis}
\big|(G_t\ast \vp_\ell\ast \psi_j)(y)\big|\lesssim \|\bar{\vp}\ast \psi_{j-\ell}\|_\infty 2^{-\ell \rho_d}\ck^{(d)}_{t,\ell}(y)\, .
\end{equation*}
In addition for the same reasons as for \eqref{bound-vp-j-ast-psi} and owing to the fact that $\psi\in \cd_L$, we have
$$\|\bar{\vp}\ast \psi_{j-\ell}\|_\infty \lesssim 2^{-L(j-\ell)}\, .$$
Following the arguments in \eqref{simp-proof-1}-\eqref{e:I-l-Lp'}, we derive
\begin{align}
2^{\ell(\al+\rho_d)}\|I\! I\! I_\ell\|_{L^p_\mu}&\lesssim \sum_{j> \ell} 2^{\ell(\al+\rho_d)}2^{-L(j-\ell)}2^{-\ell \rho_d}2^{-j\al}\Big[ 2^{j\al}\big\|\ck^{(d)}_{t,\ell} \ast |\vp_j \ast f|\big\|_{L^p_\mu}\Big]\nonumber \\
&\lesssim \sum_{ j>\ell} 2^{-(L+\al)(j-\ell)}  \Big[ 2^{j\al}\big\|\vp_j \ast f\big\|_{L^p_\mu}\Big] \lesssim \sum_{ j\geq 1} 2^{-(L+\al)j}  \Big[ 2^{(j+\ell)\al}\big\|\vp_{j+\ell} \ast f\big\|_{L^p_\mu}\Big]\, .\label{appli-rt-ii}
\end{align}
At this point, remember that the pair $(L,\vp_0)$ has been arbitrarily chosen. Therefore, we can a posteriori pick $(L,\vp_0)$ such that $L>\max(-\al,\rho_d+\al)$ (which is also consistent with the conditions exhibited in Step 1). Going back to \eqref{appli-rt-ii}, this puts us in a position to apply Jensen inequality and assert that
\begin{align*}
\sum_{\ell\ge0} 2^{\ell(\alpha+\rho_d)q}\|I\!I\!I_\ell\|^q_{L_\mu^p}&\lesssim \sum_{\ell\ge0}\bigg(\sum_{ j\geq 1} 2^{-(L+\al)j}  \Big[ 2^{(j+\ell)\al}\big\|\vp_{j+\ell} \ast f\big\|_{L^p_\mu}\Big]\bigg)^q\\
&\lesssim \sum_{\ell\ge0}\sum_{ j\geq 1} 2^{-(L+\al)j}  \Big[ 2^{(j+\ell)\al}\big\|\vp_{j+\ell} \ast f\big\|_{L^p_\mu}\Big]^q\\
&\lesssim \sum_{\ell\ge0} \Big[ 2^{\ell\al}\big\|\vp_{\ell} \ast f\big\|_{L^p_\mu}\Big]^q\lesssim \|f\|_{\cb^{\alpha,\mu}_{p,q}}^q.
\end{align*}
Accordingly, we have obtained
\begin{equation}\label{e:est-III}
\bigg(\sum_{\ell\ge0} 2^{\ell(\alpha+\rho_d)q}\|I\!I\!I_\ell\|^q_{L_\mu^p}\bigg)^{1/q}\lesssim \|f\|_{\cb^{\alpha,\mu}_{p,q}}.
\end{equation}

\noindent
\textit{Step 4: Conclusion.} By injecting \eqref{e:est-I}, \eqref{e:est-II} and \eqref{e:est-III} into \eqref{decomp-start}-\eqref{defi-i-ii}, we deduce that for some suitable $\vp_0\in\cac_c^\infty$ satisfying Assumption \eqref{condition-l-1},
$$\big\|\ell \mapsto 2^{\ell(\al+\rho_d)}\|\vp_\ell \ast \cg_tf\|_{L^p_\mu}\big\|_{\ell^q(\mathbb{N})}\lesssim\|f\|_{\cb^{\alpha,\mu}_{p,q}}\,.$$
 We are therefore in a position to apply the definition \eqref{e:Besov-norm} of Besov norm and assert that
$$\|\cg_t f\|_{B^{\al+\rho_d,\mu}_{p,q}(\vp_0)}=\big\|\ell \mapsto 2^{\ell(\al+\rho_d)}\|\vp_\ell \ast \cg_tf\|_{L^p_\mu}\big\|_{\ell^q(\mathbb{N})} \lesssim \|f\|_{B^{\al,\mu}_{p,q}(\vp_0)}\, .$$
The desired inequality \eqref{regu-effect-wave} is now a straightforward consequence of the equivalence property stated in Proposition \ref{prop:equivalence-norms}.
\end{proof}

We close this section by combining the previous results in order to obtain an appropriate regularity result for the operator $\cg_t$.

\begin{corollary}\label{cor:g-effect} We work under the same conditions as for Proposition \ref{prop:regu-effect}. Namely, consider the wave operator $\cg_t$, a constant $   \mu_\ast   >0$,    and set
$$
\rho_d:=\begin{cases}
1 & \text{if} \ d=1\\
\frac12 & \text{if} \ d=2.
\end{cases}
$$   
Let $\alpha, \mu, p, q$ be parameters such that \eqref{e:4-parameters} is fulfilled. Let $\kappa$ be an additional coefficient in $[0,   \rho_d  ]$. Then it holds that for all $0\leq s<t\leq 1$,
\begin{equation}\label{interpol}
\|\{\cg_t-\cg_s\} f\|_{\cb^{\al+\ka,\mu}_{p,q}} \lesssim |t-s|^{   \rho_d  -\ka} \|f\|_{\cb^{\al,\mu}_{p,q}}\, ,
\end{equation}
for some proportional constant that depends only on $  \mu_\ast  $.
\end{corollary}

\begin{proof}
The proof of this result, which is a standard application of interpolation properties, is included here for the sake of completeness. Indeed, write $\alpha+\kappa=(1-  \frac{\kappa}{\rho_d}  )\alpha+  \frac{\kappa}{\rho_d}  (\alpha+  \rho_d  )$.  Then apply inequality \eqref{e:interpol} with $s_0=\alpha, s_1=\alpha+   \rho_d  , \theta=  \frac{\kappa}{\rho_d}  , p_0=p_1=p, q_0=q_1=q.$ We obtain
\begin{equation}\label{e:G-interpol}
\|(\cg_t-\cg_s)f\|_{\cb^{\alpha+\kappa, \mu}_{p,q}}\le 
\|(\cg_t-\cg_s)f\|_{\cb^{\alpha, \mu}_{p,q}}^{1-  \frac{\kappa}{\rho_d}  }
\|(\cg_t-\cg_s)f\|_{\cb^{\alpha+\rho_d, \mu}_{p,q}}^{  \frac{\kappa}{\rho_d}  }\,.
\end{equation}
Next, for the first term on the right-hand side of \eqref{e:G-interpol} we apply \eqref{lipschitz}    (for $d=1$) or \eqref{e:G-Lip2} (for $d=2$)  , while for the second term we invoke \eqref{regu-effect-wave} to write 
\[\|(\cg_t-\cg_s)f\|_{\cb_{p,q}^{\alpha+   \rho_d  , \mu}}\le 
\|\cg_tf\|_{\cb_{p,q}^{\alpha+   \rho_d  , \mu}}+
\|\cg_sf\|_{\cb_{p,q}^{\alpha+   \rho_d  , \mu}}\lesssim \|f\|_{\cb^{\alpha, \mu}_{p,q}}\,.
\]
Going back to \eqref{e:G-interpol}, our conclusion easily follows. 
\end{proof}

\section{Young wave equation}\label{sec:young-eq}

This section aims to develop the Young-type integration theory  that  will allow us to  interpret and solve \eqref{e:wave} pathwisely in Section \ref{sec:wave}. For this goal, we will introduce  an  appropriate space of processes and define the Young wave integral in Section \ref{sec:Young-integral}.    Then we will    prove the wellposedness of the corresponding Young wave equation in Section \ref{sec:Young-wave}.

\subsection{Young wave integral}\label{sec:Young-integral}
 This subsection is devoted to introducing a convenient space of processes  $u$, allowing the definition of noisy integrals weighted by the wave operator. We start by giving a new piece of notation.

\begin{notation}\label{notation:mu-t} 
Throughout the section we fix two parameters $a\ge 0$ and $b>0$. We then set $\mu_t=a+bt$. 
\end{notation}

Next we introduce the set of space-time functions which will allow a proper integration theory in our noisy context. 
\begin{definition}\label{Def:space-u} Fix $T>0$ and let $(\mu_t)_{t\in [0,T]}$ be a time-dependent weight as in Notation \ref{notation:mu-t}. Consider two parameters $\gamma, \kappa\in[0,1]$. Then we define $\ce^{\ga,\ka}_{2,\infty}(T)$ as the set of functions $u:[0,T]\times \R \to \R$ for which  the following norm is finite:
\begin{equation}\label{e:norm-E}
\|u\|_{\ce^{\ga,\ka}_{2,\infty}(T)}:=\sup_{s\in [0,T]}\|u_s\|_{\cb^{\ka,\mu_s}_{2,\infty}}+\sup_{0\leq s<t \leq T} \frac{\|u_t-u_s\|_{\cb^{\ka,\mu_t}_{2,\infty}}}{|t-s|^\ga} \, .
\end{equation}
\end{definition}

For processes $u$ in a space of the form $\ce^{\ga,\ka}_{2,\infty}$, let us define some Riemann sum approximation of the integral $\int_0^t \int_{\R} G_{t-s} (x,y) u_s(y) \dot W(ds,dy).$ This definition is given below. 

\begin{definition} \label{def:Riemann-sum}
Let $P(x)=(1+|x|^2)^{-1}$ be the polynomial given in \eqref{e:weights} for $d=1$. We consider two regularity parameters $\alpha, \theta\in(0,1)$ and an integrability parameter $p>2.$ Let $\dot W$ be a noisy input which is $\theta$-H\"older continuous in time with values in the weighted space $\cb_{p,\infty}^{-\alpha, P}$; that is $\dot{W}\in \cac^\theta([0,T];\cb^{-\al,P}_{p,\infty}).$ In addition, let $\Pi^n$ be the regular dyadic partition of $[0,T]$, whose generic element is $t_m^n= mT/2^n$ for all $0\leq m\leq 2^n$. For $u$ regular enough we define the following Riemann sum based on $\Pi^n$, 
\begin{equation}\label{e:Riemann-sum}
\cj^{(n)}_{t}:=\sum_{k=0}^{m-1} \cg_{t-t^n_k}(u_{t^n_k}\, \der \dot{W}_{t^n_k t^n_{k+1}}), ~\mbox{ for } t\in (t^n_{m-1},t^n_m]\,,
\end{equation}
where we recall that $\der \dot{W}_{t^n_k t^n_{k+1}}$ stands for $\dot{W}_{t^n_{k+1}}-\dot{W}_{t^n_k}$.
\end{definition}

We now state a proposition giving natural conditions such that the Riemann sums~\eqref{e:Riemann-sum} converge in the Young sense. This result will be the basic brick towards a proper wellposedness of equation~\eqref{e:wave}.    For this statement we set, just as in Corollary \ref{cor:g-effect}, 
\begin{equation}\label{e:rho-d}
\rho_d:=\begin{cases}
1 & \text{if} \ d=1\\
\frac12 & \text{if} \ d=2.
\end{cases}
\end{equation}

\begin{proposition}\label{prop:Riemann-sum} 
Fix    two times $0<T\leq T_0$  . Let $\gamma, \theta$ be two time regularity parameters and $\kappa, \alpha$ be two space regularity parameters. Also consider $p>   d+1  $. We assume the following conditions. 
\begin{enumerate}[wide, labelwidth=!, labelindent=0pt, label=\textnormal{(\roman*)}]
\setlength\itemsep{.02in}

\item The coefficients $\gamma, \theta, \kappa, \alpha$ all sit in the interval $[0,1]$, and we have 
\begin{equation}\label{cond-ka-al-ga-the}
\ka+\al+\ga   +(1-\rho_d)   <\theta\, , \quad \ga+\theta >1\, ,\quad \ka>\al+\frac{   d  }{p}\, , \quad \ga<1-\frac{   d+1  }{p}  \, . 
\end{equation}

\item The process $u$ is an element of $ \ce^{\ga,\ka}_{2,\infty}(T)$ as introduced in Definition \ref{Def:space-u}, and $\dot W$ belongs to $\cac^\theta([0,T];\cb^{-\al,P}_{p,\infty})$.
\end{enumerate}
Let now  $\{\cj^{(n)}; n\ge 1\}$ be the sequence defined by \eqref{e:Riemann-sum}. Then we have that $\cj^{(n)}$ converges in $\ce_{2,\infty}^{\gamma, \kappa}(T)$.   We denote 
\begin{equation}\label{e:Young-integral}
\lim_{n\to\infty} \cj_t^{(n)} =: \int_0^t \cg_{t-r}(u_r \, d\dot W_r)\, .
\end{equation}
Moreover,  the Young integral~\eqref{e:Young-integral} verifies 
\begin{equation}
\bigg\|\int_0^. \cg_{.-r}(u_r \, d\dot W_r)\bigg\|_{\ce^{\ga,\ka}_{2,\infty}(T)}\leq  \big\|\cg_.(u_0\, (\dot{W}_{T}-\dot{W}_0))\big\|_{\ce_{2,\infty}^{\gamma,\kappa}(T)} 
+   c_{T_0}  \|\dot W\|_{\cac^\theta([0,T];\cb^{-\al,P}_{p,\infty})}  \|u\|_{\ce^{\gamma,\kappa}_{2,\infty}(T)} \, ,\label{boun-int-you-wav}
\end{equation}
 where $c_{T_0}>0$ does not depend on $T$, $u$ and $\dot{W}$.
\end{proposition}

\begin{remark} As the reader will see, we will in fact prove a slightly stronger version of relation~\eqref{boun-int-you-wav}:  there exists a \textit{finite} constant $q=q(\ga,\ka,\theta,\al,p)\geq 1$ such that \begin{equation}\label{boun-int-you-wav'}
\begin{split}
&\bigg\|\int_0^. \cg_{.-r}(u_r \, d\dot W_r)\bigg\|_{\ce^{\ga,\ka}_{2,\infty}(T)}\\
&\leq  \big\|\cg_.(u_0\, (\dot{W}_{T}-\dot{W}_0))\big\|_{\ce_{2,\infty}^{\gamma,\kappa}(T)} + c_{T_0}\|\dot W\|_{\cac^\theta([0,T];\cb^{-\al,P}_{p,\infty})} \Big(\int_0^T \|u\|_{\ce^{\ga,\ka}_{2,\infty}(r)}^q \, dr\Big)^{\frac{1}{q}}  \, ,
\end{split}
\end{equation}
 where $c_{T_0}>0$ does not depend on $T$, $u$ and $\dot{W}$. Inequality \eqref{boun-int-you-wav'} will be crucial in our fixed point argument for the existence-uniqueness of  a solution to  \eqref{e:wave}. 
\end{remark}

\begin{proof}[Proof of Proposition \ref{prop:Riemann-sum}]
We shall prove the convergence of $\cj^{(n)}$ thanks to an upper bound on the difference $\cj^{(n+1)} -\cj^{(n)}$. To this aim, let us denote $\delta f_{st}$ for the increment $f_t-f_s$ of a function $f$. Furthermore, for the sake of clarity, consider two dyadic points $s=t^n_\ell$ and $t=t^n_m$ for $\ell \le m$. Then an elementary manipulation on the expression \eqref{e:Riemann-sum} for $\cj^{(n)}$ reveals that 
\begin{equation}\label{e:decom-J}
\delta\cj_{st}^{(n+1)}-\delta\cj_{st}^{(n)} =I^{(n)}_{st}+I\! I^{(n)}_{st}+I\! I\! I^{(n)}_{st}+I\! V^{(n)}_{st}\, ,
\end{equation}
where the terms $I^{(n)}_{st}, \dots, I\! V^{(n)}_{st}$ are defined by (writing $t_{j}=t_{j}^{n+1}$ for notational sake), 
\begin{align*}
I^{(n)}_{st}&= \sum_{k=\ell}^{m-1} \cg_{t-t_{2k}}\big( \der u_{t_{2k}t_{2k+1}} \, \der \dot{W}_{t_{2k+1}t_{2k+2}}\big);\\
I\!I^{(n)}_{st}&=\sum_{k=\ell}^{m-1} \big\{\cg_{t-t_{2k+1}}-\cg_{t-t_{2k}}\big\} \big(u_{t_{2k+1}} \, \der \dot{W}_{t_{2k+1}t_{2k+2}}\big);\\
I\!I\!I_{st}^{(n)}&=\sum_{k=0}^{\ell-1} \big\{\cg_{t-t_{2k}}-\cg_{s-t_{2k}}\big\}\big( \der u_{t_{2k}t_{2k+1}} \, \der \dot{W}_{t_{2k+1}t_{2k+2}}\big);\\
I\!V_{st}^{(n)}&= \sum_{k=0}^{\ell-1} \big\{\cg_{t-t_{2k+1}}-\cg_{t-t_{2k}}-\cg_{s-t_{2k+1}}+\cg_{s-t_{2k}}\big\} \big(u_{t_{2k+1}} \, \der \dot{W}_{t_{2k+1}t_{2k+2}}\big).
\end{align*}
We now upper bound the four terms above.

\noindent
\textit{Step 0: A general bound on products.}   Recall that the (fixed) parameters $a$ and $b$ have been introduced in Notation \ref{notation:mu-t}. Setting
$$\mu_\ast=\mu_\ast(T_0,p):=\max\bigg(\frac{a+bT_0}{1-\frac{2}{p}},\frac{pbT_0}{2} \bigg),$$
it is readily checked that for all $0\leq  s<t\leq T$,
\begin{equation}\label{def-mu-ast}
\mu_s \leq \mu_\ast \Big(1-\frac{2}{p}\Big) \quad \text{and}\quad \frac{p}{2}(\mu_t-\mu_s)\leq \mu_\ast. 
\end{equation}

Consider our noisy input $\dot W$ and a function $f\in \cb_{2, \infty}^{\kappa, \mu_s}$ for $0\le s<t<T$. Then we claim that for $s\le r < v \le t$ we have 
\begin{align}
\big\| f\cdot \delta \dot{W}_{rv}\big\|_{\cb^{-\al,\mu_t}_{2,\infty}} &\lesssim \big\| f \big\|_{\cb^{\ka,\mu_{s}}_{2,\infty}}\frac{|v-r|^{\theta}}{|t-s|^{\frac{ d+1}{p}}}\big\|  \dot{W}\big\|_{\cac^\theta([0,T];\cb^{-\al,P}_{p,\infty})}\, ,\label{estim-prod-f-w}
\end{align}
for some proportional constant that depends only on $\mu_\ast$, that is only on $T_0$ and $p$ (note that we will no longer explicitly indicate such a dependence on $(T_0,p)$ for the constants arising in the rest of the proof).

%\smallskip

 To prove \eqref{estim-prod-f-w}, observe first that, owing to \eqref{def-mu-ast}, we are in a position to apply our product rule \eqref{e:product-rule} and assert that
\[\big\| f\cdot \delta \dot{W}_{rv}\big\|_{\cb^{-\al,\mu_t}_{2,\infty}}\lesssim \big\| f \big\|_{\cb^{\ka,\mu_{s}}_{2,\infty}}\big\| \delta \dot{W}_{rv}\big\|_{\cb^{-\al,\varepsilon_{st}}_{p,\infty}}\,,
\]
where we have set $\varepsilon_{st}=\frac p2 (\mu_t-\mu_s)\leq \mu_\ast$. Next we invoke Lemma~
\ref{lem:change-weight} for the right-hand side above, which allows to write 
\[ \big\| f\cdot \delta \dot{W}_{rv}\big\|_{\cb^{-\al,\mu_t}_{2,\infty}}\lesssim \big\| f \big\|_{\cb^{\ka,\mu_{s}}_{2,\infty}}\frac{1}{(\mu_t-\mu_s)^{\frac{ d+1}{p}}}\big\| \delta \dot{W}_{rv}\big\|_{\cb^{-\al,P}_{p,\infty}}
\,,\]
from which \eqref{estim-prod-f-w} is easily deduced.

\noindent
\textit{Step 1: Bound for $I_{st}^{(n)}$.} For the term $I_{st}^{(n)}$ in the right-hand side of \eqref{e:decom-J},  let us first simply write
\begin{equation}\label{e:I-1}
 \big\|I^{(n)}_{st}\big\|_{\cb^{\ka,\mu_t}_{2,\infty}}\leq \sum_{k=\ell}^{m-1} \big\|\cg_{t-t_{2k}}\big( \der u_{t_{2k}t_{2k+1}} \, \der \dot{W}_{t_{2k+1}t_{2k+2}}\big)\big\|_{\cb^{\ka,\mu_t}_{2,\infty}}\,.
 \end{equation}
Then we apply \eqref{interpol}, interpolating from a regularity $(-\alpha)$ to a regularity $\kappa$ and considering the times $t:=t-t_{2k}$ and $s:=0$ (recall that $G_{0}=0$ owing to~\eqref{e:G}). This yields
\begin{equation}\label{e:I-2} \big\|I^{(n)}_{st}\big\|_{\cb^{\ka,\mu_t}_{2,\infty}}\lesssim \sum_{k=\ell}^{m-1}|t-t_{2k}|^{ \rho_d-(\ka+\al)} \big\| \der u_{t_{2k}t_{2k+1}}\, \der \dot{W}_{t_{2k+1}t_{2k+2}}\big\|_{\cb^{-\al,\mu_t}_{2,\infty}}\,. 
\end{equation}
We are now in a position to resort to our general estimate \eqref{estim-prod-f-w}, which gives
\begin{equation}\label{e:I-3}
\big\|I^{(n)}_{st}\big\|_{\cb^{\ka,\mu_t}_{2,\infty}}\lesssim  \sum_{k=\ell}^{m-1}|t-t_{2k+1}|^{ \rho_d-(\ka+\al)}\big\|\der u_{t_{2k}t_{2k+1}} \big\|_{\cb^{\ka,\mu_{t_{2k+1}}}_{2,\infty}}\frac{2^{-n\theta}}{|t-t_{2k+1}|^{\frac{ d+1}{p}}} \,
\|  \dot{W}\|_{\cac^\theta([0,T];\cb^{-\al,P}_{p,\infty})}\,.
\end{equation}
Invoking the Definition \eqref{e:norm-E} for the norm in $\ce_{2,\infty}^{\gamma, \mu}$ and rearranging terms, we thus get 
\begin{equation}\label{e:I-4}
\big\|I^{(n)}_{st}\big\|_{\cb^{\ka,\mu_t}_{2,\infty}}\lesssim  2^{-n(\ga+\theta-1)} \|\dot{W}\|_{\cac^\theta([0,T];\cb^{-\al,P}_{p,\infty})} \bigg(2^{-n} \sum_{k=\ell}^{m-1}|t-t_{2k+1}|^{ \rho_d-(\ka+\al)-\frac{ d+1}{p}}\|u\|_{\ce^{\ga,\ka}_{2,\infty}(t_{2k+1})}\bigg)\,.
\end{equation}
The last expression in the right-hand side above can be upper bounded by the continuous integral 
\begin{equation}\label{e:bd-int1}
\int_s^t (t-v)^{\rho_d-(\kappa+\alpha)-\frac{ d+1}{p}}\|u\|_{\ce^{\ga,\ka}_{2,\infty}(v)} \, dv.
\end{equation}
Therefore \eqref{e:I-4} can be recast as 
\begin{equation}\label{e:bd-I}
 \big\|I^{(n)}_{st}\big\|_{\cb^{\ka,\mu_t}_{2,\infty}}
 \lesssim  
 2^{-n(\gamma+\theta-1)} \|\dot{W}\|_{\cac^\theta([0,T];\cb^{-\al,P}_{p,\infty})}\, 
 \int_s^t (t-v)^{\rho_d-(\kappa+\alpha)-\frac{ d+1 }{p}}\|u\|_{\ce^{\ga,\ka}_{2,\infty}(v)} \, dv\, .
 \end{equation}
We now wish to apply H\"older's inequality to the integral term \eqref{e:bd-int1} above. Namely, for $r_1, q_1$ such that $\frac{1}{r_1}+\frac{1}{q_1}=1$, we write
\begin{multline} \label{e:int-holder}
\int_s^t (t-v)^{\rho_d-(\kappa+\alpha)-\frac{ d+1}{p}}\|u\|_{\ce^{\ga,\ka}_{2,\infty}(v)} \, dv\, \\
\le \bigg(\int_s^t (t-v)^{r_1(\rho_d -(\kappa+\alpha)-\frac{ d+1}{p})}\, dv\bigg)^{\frac{1}{r_1}} \bigg( \int_s^t \|u\|_{\ce^{\ga,\ka}_{2,\infty}(v)}^{q_1} \, dv\bigg)^{\frac{1}{q_1}}\,. 
 \end{multline}
 In order to make the right-hand side of \eqref{e:int-holder} finite, we need to have $r_1(\rho_d-(\kappa+\alpha)-\frac{d+1}{p})>-1$. Now observe that thanks to the conditions in \eqref{cond-ka-al-ga-the} we can successively guarantee that
$$\ga- \rho_d+\ka+\al+\frac{ d+1 }{p}<\theta-\lp 1-\frac{ d+1 }{p}\rp<\theta-\ga\leq 1.$$
As a result, we can pick $r_1>1$ such that
\begin{equation}\label{defi-par-r}
\ga- \rho_d+\ka+\al+\frac{ d+1}{p}<\frac{1}{r_1}
\quad\Longrightarrow\quad
r_1\lp   \rho_d  -(\kappa+\alpha)-\frac{   d+1  }{p}\rp>-1+ r_{1}\ga.
\end{equation}
With this set of parameters in hand, inequality \eqref{e:int-holder} reads
\begin{eqnarray*}
 \int_s^t (t-v)^{  \rho_d  -(\kappa+\alpha)-\frac{   d+1  }{p}}\|u\|_{\ce^{\ga,\ka}_{2,\infty}(v)} \, dv
&\lesssim  &
|t-s|^{  \rho_d  -(\kappa+\alpha)-\frac{   d+1  }{p}+\frac{1}{r_1}}\bigg( \int_0^T \|u\|_{\ce^{\ga,\ka}_{2,\infty}(v)}^{q_1} \, dv\bigg)^{\frac{1}{q_1}}\\
&\lesssim &    |t-s|^{\gamma}  \bigg( \int_0^T \|u\|_{\ce^{\ga,\ka}_{2,\infty}(v)}^{q_1} \, dv\bigg)^{\frac{1}{q_1}} \, ,
\end{eqnarray*}
where we have invoked \eqref{defi-par-r} to derive the second step above.  Plugging this inequality into \eqref{e:bd-I} and setting $\varepsilon_1:=\gamma+\theta-1>0$ (thanks to \eqref{cond-ka-al-ga-the}), we end up with 
 \begin{equation}\label{a2}
\big\|I^{(n)}_{st}\big\|_{\cb^{\ka,\mu_t}_{2,\infty}}\lesssim  
   2^{-n\varepsilon_1} |t-s|^{\ga}    \|\dot{W}\|_{\cac^\theta([0,T];\cb^{-\al,P}_{p,\infty})}\, 
\bigg( \int_0^T \|u\|_{\ce^{\ga,\ka}_{2,\infty}(v)}^{q_1} \, dv\bigg)^{\frac{1}{q_1}}\,.
\end{equation}

\noindent
\textit{Step 2: Bound for $I\!I_{st}^{(n)}$.}  The computations are similar to those of Step 1, and we will just highlight the main differences. First the equivalent of \eqref{e:I-1} is 
\[\big\|I\!I^{(n)}_{st}\big\|_{\cb^{\ka,\mu_t}_{2,\infty}}\leq \sum_{k=\ell}^{m-1} \big\|\big\{\cg_{t-t_{2k+1}}-\cg_{t-t_{2k}}\big\} \big(u_{t_{2k+1}} \, \der \dot{W}_{t_{2k+1}t_{2k+2}}\big)\big\|_{\cb^{\ka,\mu_t}_{2,\infty}}\,.\]
Then we invoke the regularity property \eqref{interpol} and the fact that $t_{2k+1}-t_{2k}=2^{-(n+1)}$ in order to get 
\[   \big\|I\! I^{(n)}_{st}\big\|_{\cb^{\ka,\mu_t}_{2,\infty}}    \lesssim  2^{-n(   \rho_d  -(\ka+\al))} \sum_{k=\ell}^{m-1} \big\|u_{t_{2k+1}} \, \der \dot{W}_{t_{2k+1}t_{2k+2}}\big\|_{\cb^{-\al,\mu_t}_{2,\infty}}\,,\]
which is parallel to \eqref{e:I-2}. We now repeat the arguments in \eqref{e:I-3}-\eqref{e:I-4}, which yields
\begin{multline}\label{e:II-1}
   \big\|I\! I^{(n)}_{st}\big\|_{\cb^{\ka,\mu_t}_{2,\infty}}    \\
   \lesssim  2^{-n(\theta   +\rho_d-1  -(\kappa+\alpha))}  \|\dot{W}\|_{\cac^\theta([0,T];\cb^{-\al,P}_{p,\infty})}  \bigg( 2^{-n} \sum_{k=\ell}^{m-1}\frac{1}{|t-t_{2k+1}|^{\frac{   d+1  }{p}}}\|u\|_{\ce^{\ga,\ka}_{2,\infty}(t_{2k+1})}\bigg)\,.
\end{multline}
The Riemann type sum in \eqref{e:II-1} can be upper bounded by $\int_s^t (t-v)^{-\frac{   d+1  }{p}}\|u\|_{\ce^{\ga,\ka}_{2,\infty}(v)} \, dv$. Besides, recall from \eqref{cond-ka-al-ga-the} that we have assumed $1-\frac{   d+1  }{p}>\gamma$, and so there exists $r_2>1$ such that $\frac{1}{r_2}-\frac{   d+1  }{p}>\ga$. Denoting by $q_2$ the (finite) conjugate of $r_2$, we can go back to \eqref{e:II-1} and deduce
\begin{align*}%\label{e:II-2}
&   \big\|I\! I^{(n)}_{st}\big\|_{\cb^{\ka,\mu_t}_{2,\infty}}    \lesssim  2^{-n(\theta   +\rho_d-1  -(\kappa+\alpha))} \|\dot{W}\|_{\cac^\theta([0,T];\cb^{-\al,P}_{p,\infty})} \int_s^t (t-v)^{-\frac{   d+1  }{2}}\|u\|_{\ce^{\ga,\ka}_{2,\infty}(v)} \, dv\\
&\lesssim  2^{-n(\theta   +\rho_d-1  -(\kappa+\alpha))} \|\dot{W}\|_{\cac^\theta([0,T];\cb^{-\al,P}_{p,\infty})}  \bigg(\int_s^t (t-v)^{-\frac{r_2   (d+1)  }{p}}\, dv\bigg)^{\frac{1}{r_2}} \bigg( \int_s^t \|u\|_{\ce^{\ga,\ka}_{2,\infty}(v)}^{q_2} \, dv\bigg)^{\frac{1}{q_2}}\\
&\lesssim  2^{-n(\theta   +\rho_d-1  -(\kappa+\alpha))} \|\dot{W}\|_{\cac^\theta([0,T];\cb^{-\al,P}_{p,\infty})}  |t-s|^{\frac{1}{r_2}-\frac{   d+1  }{p}}\bigg( \int_0^T \|u\|_{\ce^{\ga,\ka}_{2,\infty}(v)}^{q_2} \, dv\bigg)^{\frac{1}{q_2}},
\end{align*}
which finally yields
\begin{equation}\label{e:II-3}
   \big\|I\! I^{(n)}_{st}\big\|_{\cb^{\ka,\mu_t}_{2,\infty}}    \lesssim  2^{-n\varepsilon_2} \|u\|_{\ce^{\ga,\ka}_{2,\infty}} \|\dot{W}\|_{\cac^\theta([0,T];\cb^{-\al,P}_{p,\infty})}    |t-s|^\gamma  \bigg( \int_0^T \|u\|_{\ce^{\ga,\ka}_{2,\infty}(v)}^{q_2} \, dv\bigg)^{\frac{1}{q_2}},
\end{equation}
where $\varepsilon_2:=\theta   +\rho_d-1  -(\kappa+\alpha)>\gamma\geq 0$ (thanks to \eqref{cond-ka-al-ga-the}).

\noindent
\textit{Step 3: Bound for $I\!I\!I_{st}^{(n)}$.} The term $I\!I\!I_{st}^{(n)}$ in the right-hand side of \eqref{e:decom-J} is treated again as in Step 1 and Step 2. Namely applying once again \eqref{interpol} we get 
\[\big\|I\!I\! I^{(n)}_{st}\big\|_{\cb^{\ka,\mu_t}_{2,\infty}}\lesssim  |t-s|^{   \rho_d  -(\ka+\al)}\sum_{k=0}^{\ell-1} \big\| \der u_{t_{2k}t_{2k+1}} \, \der \dot{W}_{t_{2k+1}t_{2k+2}}\big\|_{\cb^{-\al,\mu_t}_{2,\infty}}\,.
\]
Then we invoke the fact that $\mu_t<\mu_{t_{2k+1}}$ for all $k<\ell-1$ and we take advantage of the time regularity of $u$ and $\dot W$, mimicking again \eqref{e:I-3}-\eqref{e:I-4}. This allows to write
\begin{multline*}
\big\|I\!I\! I^{(n)}_{st}\big\|_{\cb^{\ka,\mu_t}_{2,\infty}} \\
\lesssim 
2^{-n(\theta+\gamma-1)} \|\dot{W}\|_{\cac^\theta([0,T];\cb^{-\al,P}_{p,\infty})}|t-s|^{   \rho_d  -(\kappa+\alpha)} \bigg( 2^{-n}   \sum_{k=0}^{\ell-1}\frac{1}{|t-t_{2k+1}|^{\frac{   d+1  }{p}}}\|u\|_{\ce^{\ga,\ka}_{2,\infty}(t_{2k+1})}\bigg)\,.
\end{multline*}
Proceeding as in Step 2, we can upper bound the quantity into brackets by
\begin{align*}
\int_{0}^{t}(t-r)^{-\frac{   d+1  }{2}}\|u\|_{\ce^{\ga,\ka}_{2,\infty}(r)} dr&\leq T^{   \frac{1}{r_2}-\frac{d+1}{p}  } \bigg( \int_0^T \|u\|_{\ce^{\ga,\ka}_{2,\infty}(v)}^{q_3} \, dv\bigg)^{\frac{1}{q_3}}
\end{align*}
where $q_3=q_2$.
We end up with 
\begin{equation}\label{e:III-1}
 \big\|I\!I\! I^{(n)}_{st}\big\|_{\cb^{\ka,\mu_t}_{2,\infty}} \lesssim 2^{-n\varepsilon_3}|t-s|^\gamma \|\dot{W}\|_{\cac^\theta([0,T];\cb^{-\al,P}_{p,\infty})}\bigg( \int_0^T \|u\|_{\ce^{\ga,\ka}_{2,\infty}(v)}^{q_3} \, dv\bigg)^{\frac{1}{q_3}}\,,
 \end{equation}
where $\varepsilon_3:=\gamma+\theta-1>0$.

\noindent
\textit{Step 4: Bound for $I\!V_{st}^{(n)}$.}  For the term $I\!V_{st}^{(n)}$, we will bound the increment of $\cg$ in two different ways. Specifically, alleviate notations by setting 
\[\mathcal V_{st}^k:= \big\{\cg_{t-t_{2k+1}}-\cg_{t-t_{2k}}-\cg_{s-t_{2k+1}}+\cg_{s-t_{2k}}\big\} \big(u_{t_{2k+1}} \, \der \dot{W}_{t_{2k+1}t_{2k+2}}\big)\,.\]
Then on the one hand we have 
\begin{align}\label{e:V-1}
\big\|\mathcal V^{k}_{st}\big\|_{\cb^{\ka,\mu_t}_{2,\infty}}&\leq   \big\|\big\{\cg_{t-t_{2k+1}}-\cg_{s-t_{2k+1}}\big\} \big(u_{t_{2k+1}} \, \der \dot{W}_{t_{2k+1}t_{2k+2}}\big)\big\|_{\cb^{\ka,\mu_t}_{2,\infty}}\notag\\
&\qquad  + \big\|\big\{\cg_{t-t_{2k}}-\cg_{s-t_{2k}}\big\} \big(u_{t_{2k+1}} \, \der \dot{W}_{t_{2k+1}t_{2k+2}}\big)\big\|_{\cb^{\ka,\mu_t}_{2,\infty}}\notag\\
&\lesssim |t-s|^{   \rho_d  -(\kappa+\alpha)} \|u_{t_{2k+1}} \, \der \dot{W}_{t_{2k+1}t_{2k+2}}\big\|_{\cb^{-\alpha,\mu_t}_{2,\infty}}\,,
\end{align}
where we have invoked \eqref{interpol} for the second inequality. On the other hand, still using \eqref{interpol} we also have 
\begin{align}\label{e:V-2}
\big\|\mathcal V^{k}_{st}\big\|_{\cb^{\ka,\mu_t}_{2,\infty}}&\leq   \big\|\big\{\cg_{t-t_{2k+1}}-\cg_{t-t_{2k}}\big\} \big(u_{t_{2k+1}} \, \der \dot{W}_{t_{2k+1}t_{2k+2}}\big)\big\|_{\cb^{\ka,\mu_t}_{2,\infty}}\notag\\
&\qquad  + \big\|\big\{\cg_{s-t_{2k+1}}-\cg_{s-t_{2k}}\big\} \big(u_{t_{2k+1}} \, \der \dot{W}_{t_{2k+1}t_{2k+2}}\big)\big\|_{\cb^{\ka,\mu_t}_{2,\infty}}\notag\\
&\lesssim 2^{-n(   \rho_d  -(\kappa+\alpha))} \|u_{t_{2k+1}} \, \der \dot{W}_{t_{2k+1}t_{2k+2}}\big\|_{\cb^{-\alpha,\mu_t}_{2,\infty}}\,.
\end{align}
We now introduce an additional parameter $\beta\in(0,1)$, whose exact value will be specified later on. Combining \eqref{e:V-1} and \eqref{e:V-2}, we get 
\[\big\|I\! V^{(n)}_{st}\big\|_{\cb^{\ka,\mu_t}_{2,\infty}}\lesssim \sum_{k=0}^{\ell-1} \big[|t-s|^{   \rho_d  -(\ka+\al)} \big]^\beta  \big[2^{-n(  \rho_d  -(\ka+\al))} \big]^{1-\beta}\big\|u_{t_{2k+1}} \, \der \dot{W}_{t_{2k+1}t_{2k+2}}\big\|_{\cb^{-\al,\mu_t}_{2,\infty}}\,.
\]
Along the same lines as for \eqref{e:I-3}-\eqref{e:I-4}, we then let the reader check that we have 
\begin{multline}\label{e:IV-1}
\big\|I\! V^{(n)}_{st}\big\|_{\cb^{\ka,\mu_t}_{2,\infty}}\lesssim 2^{-n(   \rho_d  -(\kappa+\alpha))(1-\beta)}|t-s|^{(   \rho_d  -(\kappa+\alpha))\beta} \\
\times
 \|\dot{W}\|_{\cac^\theta([0,T];\cb^{-\al,P}_{p,\infty})}    2^{-n(\theta-1)} \bigg(2^{-n}  \sum_{k=0}^{\ell-1} \frac1{|t-t_{2k+1}|^{\frac{   d+1  }{p}}}\|u\|_{\ce^{\ga,\ka}_{2,\infty}(t_{2k+1})}\bigg)\,.
\end{multline}
We choose $\beta\in(0,1)$ such that $(  \rho_d  -(\kappa+\alpha))\beta=\gamma$, which is compatible with our assumption~\eqref{cond-ka-al-ga-the}. Bounding the Riemann sum in \eqref{e:IV-1} just as in Step 3, that is by
$$\int_{0}^{t}(t-r)^{-\frac{   d+1  }{2}}\|u\|_{\ce^{\ga,\ka}_{2,\infty}(r)} dr\leq T^{   \frac{1}{r_2}-\frac{   d+1  }{p}  } \bigg( \int_0^T \|u\|_{\ce^{\ga,\ka}_{2,\infty}(v)}^{q_4} \, dv\bigg)^{\frac{1}{q_4}},$$
where $q_4=q_3$, we obtain 
\begin{equation}\label{e:IV-2}
\big\|I\! V^{(n)}_{st}\big\|_{\cb^{\ka,\mu_t}_{2,\infty}}\lesssim 2^{-n\varepsilon_4}|t-s|^{\gamma} \|\dot{W}\|_{\cac^\theta([0,T];\cb^{-\al,P}_{p,\infty})} \bigg( \int_0^T \|u\|_{\ce^{\ga,\ka}_{2,\infty}(v)}^{q_4} \, dv\bigg)^{\frac{1}{q_4}}\,,
\end{equation}
with the notation $\varepsilon_4=   \rho_d-(\ka+\al)-\ga+\theta-1  $. Notice that $\varepsilon_4>0$ thanks to \eqref{cond-ka-al-ga-the}.

\noindent\textit{Step 5: Conclusion.}  Recall the decomposition \eqref{e:decom-J}. Then owing to \eqref{a2}, \eqref{e:II-3}, \eqref{e:III-1} and~\eqref{e:IV-2}, we have obtained that for all $0\le s<t \le T$,
\begin{equation} \label{bou-step-5}
\Big\|\delta \cj_{st}^{(n+1)}-\delta\cj_{st}^{(n)}\Big\|_{\cb^{\ka,\mu_t}_{2,\infty}}\lesssim    2^{-n\varepsilon} |t-s|^\ga\|\dot{W}\|_{\cac^\theta([0,T];\cb^{-\al,P}_{p,\infty})}   \bigg( \int_0^T \|u\|_{\ce^{\ga,\ka}_{2,\infty}(v)}^{q} \, dv\bigg)^{\frac{1}{q}}\, ,
\end{equation}
where $\varepsilon:=\min_{i=1,\ldots,4}\varepsilon_i>0$ and $q:=\max_{i=1,\ldots,4} q_i$. Using the fact that $\cj^{(n)}_0=0$, the bound \eqref{bou-step-5} also allows us to assert that for all $0\leq s\leq T$,
$$
\Big\|\cj_{s}^{(n+1)}-\cj_{s}^{(n)}\Big\|_{\cb^{\ka,\mu_s}_{2,\infty}}\lesssim 2^{-n\varepsilon} \|\dot{W}\|_{\cac^\theta([0,T];\cb^{-\al,P}_{p,\infty})} \bigg( \int_0^T \|u\|_{\ce^{\ga,\ka}_{2,\infty}(v)}^{q} \, dv\bigg)^{\frac{1}{q}}\, ,
$$
 Recalling the definition \eqref{e:norm-E} of the norm $\|\cdot\|_{\ce_{2,\infty}^{\gamma,\kappa}(T)}$, we  have thus established that
\begin{equation}\label{toward-cauchy}
\|\cj^{(n+1)}-\cj^{(n)}\|_{\ce_{2,\infty}^{\gamma,\kappa}(T)}\lesssim 2^{-n\varepsilon}\|\dot{W}\|_{\cac^\theta([0,T];\cb^{-\al,P}_{p,\infty})}\bigg( \int_0^T \|u\|_{\ce^{\ga,\ka}_{2,\infty}(v)}^{q} \, dv\bigg)^{\frac{1}{q}}\,.
\end{equation}
This shows that $\cj^{(n)}$ is a Cauchy sequence in $\ce^{\ga,\ka}_{2,\infty}(T)$, which concludes the proof of the convergence statement.

Moreover, denoting (temporarily) the limit of $\cj^{(n)}$ by $\cj$, we can use again \eqref{toward-cauchy} to deduce
$$
\|\cj\|_{\ce_{2,\infty}^{\gamma,\kappa}(T)}\leq \|\cj^{(0)}\|_{\ce_{2,\infty}^{\gamma,\kappa}(T)}+c_{T_0}\,\|\dot{W}\|_{\cac^\theta([0,T];\cb^{-\al,P}_{p,\infty})} \bigg( \int_0^T \|u\|_{\ce^{\ga,\ka}_{2,\infty}(v)}^{q} \, dv\bigg)^{\frac{1}{q}} \, ,
$$
for some constant $   c_{T_0}  >0$, which exactly corresponds to the desired inequality \eqref{boun-int-you-wav'}. Also remember that \eqref{boun-int-you-wav'} is stronger than \eqref{boun-int-you-wav}, and thus \eqref{boun-int-you-wav} holds as well. This finishes the proof. 
\end{proof}

\subsection{Wellposedness of the Young wave equation for $\mathbf{d=1,2}$}\label{sec:Young-wave}

Recall that the wave operator $\cg_t$ has been introduced in Section \ref{sec:wave-kernel}. With the preliminary results of Section \ref{sec:Young-integral} in hand, equation~\eqref{e:wave} with initial conditions $u_0, u_1$ will now be written in the following mild sense: 
\begin{equation}\label{equa-young-wellposed}
u_t:=(\partial_t \cg)_t u_0+\cg_tu_1+\int_0^t \cg_{t-r}(u_r \, d\dot W_r),
\end{equation}
In \eqref{equa-young-wellposed} the integral is interpreted in the Young sense, that is via Proposition \ref{prop:Riemann-sum}.

The initial conditions in \eqref{equa-young-wellposed} have to satisfy some standard smoothness conditions. Those conditions will be expressed in the usual Sobolev scale, for which we introduce a new notation. 

\begin{notation}\label{def:Hs}
For every $s\geq 0$, we denote by $\ch^s=\ch^s(\R^{   d  })$ the usual Sobolev space of order $s$, that is the set of functions $f\in L^2(\R^{   d  })$ such that
$$\|f\|_{\ch^s}^2:=\int_{\R^{   d  }}d\xi \, \{1+|\xi|^2\}^s |\hat{f}(\xi)|^2 \, < \, \infty.$$ 
\end{notation}

We now turn to our main abstract existence and uniqueness result.    Let us recall that we have set
$$
\rho_d:=\begin{cases}
1 & \text{if} \ d=1\\
\frac12 & \text{if} \ d=2.
\end{cases}
$$

\begin{theorem}\label{theo:wellposedness}
Assume $d   \in \{1,2\}  $.  Fix an arbitrary time $T>0$, as well as parameters $\gamma, \theta, \kappa, \alpha\in [0,1]$ and $p>2$ satisfying the conditions in \eqref{cond-ka-al-ga-the}. Assume that $\dot W$ belongs to $\cac^\theta([0,T];\cb^{-\al,P}_{p,\infty})$ and pick $(u_0,u_1)\in \ch^{1+\ka}\times \ch^\ka$. 
Then equation \eqref{equa-young-wellposed} admits a unique solution $u$ in the space $\ce^{\ga,\ka}_{2,\infty}(T)$ introduced in Definition \ref{Def:space-u}.
\end{theorem}

Before proving Theorem \ref{theo:wellposedness}, we start with a technical lemma giving a control on the terms related to the initial conditions. 

\begin{lemma}\label{lem:ini-condi}
Under the assumptions of Theorem \ref{theo:wellposedness} and for all $t\in [0,T]$, set
\begin{equation}\label{e:u-0}
u^{(0)}_t:=(\partial_t \cg)_t u_0+\cg_t u_1 \ .
\end{equation}
Then $u^{(0)}$ is an element of  $\ce^{\ga,\ka}_{2,\infty}(T)$ and we have
\begin{equation}\label{e:bd-u-0}
\big\|u^{(0)}\big\|_{\ce^{\ga,\ka}_{2,\infty}(T)} \lesssim \big\|u_0\big\|_{\ch^{1+\ka}}+\big\|u_1\big\|_{\ch^{\ka}}.
\end{equation}
\end{lemma}

\begin{proof}
Let $\vp\in \cd_1$ denote the function involved in our definition of the spaces $\cb^{\ka,\mu}_{p,q}$ (see Definition~\ref{def:Besov} and Notation \ref{notation:pmb-vp}). As a preliminary step, observe that owing to assumption \eqref{condition-l-1} and the fact that $\vp$ has compact support, one has 
$$|\widehat{\vp}(\xi)|=\bigg|\int_{\R^d}dx \, e^{-\imath \xi\cdot x} \vp(x)\bigg|=\bigg|\int_{\R^d}dx \, \{e^{-\imath \xi\cdot x}-1\} \vp(x)\bigg|\lesssim |\xi|^\ka,$$
where we have used the fact that $\kappa\in(0,1)$. Hence for every $j\geq 1$,
\begin{equation}\label{estim-hat-phi-j}
|\widehat{\vp_j}(\xi)|=|\widehat{\vp}(2^{-j}\xi)|\lesssim 2^{-j\ka}|\xi|^\ka .
\end{equation}
Now the weight $w=e^{-\mu|x|}$ is smaller than 1. Thus one can trivially bound the norm in $L^2_{\mu_t}$ by the usual norm in $L^2(\R^d)$. Hence invoking \eqref{e:Fourier-G} we get 
\begin{align*}
\big\|\vp_{j}\ast \{(\partial_t \cg)_t -(\partial_t \cg)_s\} u_0\big\|_{L^2_{\mu_t}} &\lesssim \big\|\big(\vp_j \ast(\{(\partial_t G)_t -(\partial_t G)_s\} \ast  u_0\big)\big\|_{L^2}\\
&\lesssim \bigg( \int_{\R^d}d\xi \, \big| \widehat{\vp_j}(\xi)  \big|^2\big|\cos(t|\xi|)-\cos(s|\xi|)\big|^2  \big|\widehat{u_0}(\xi)\big|^2\bigg)^{\frac12} \, .
\end{align*}
Combining this inequality with \eqref{estim-hat-phi-j}, this yields
\begin{eqnarray}\label{e:bd-vp-G-u0}
\big\|\vp_{j}\ast \{(\partial_t \cg)_t -(\partial_t \cg)_s\} u_0\big\|_{L^2_{\mu_t}}
&\lesssim& 
2^{-j\ka} |t-s|\bigg( \int_{\R^d}d\xi \, |\xi|^{2+2\ka} \big|\widehat{u_0}(\xi)\big|^2\bigg)^{\frac12} \notag\\
&\lesssim& 
2^{-j\ka} |t-s|^\ga \|u_0\|_{\ch^{1+\ka}}\,.
\end{eqnarray}
Plugging \eqref{e:bd-vp-G-u0} into the definition \eqref{e:Besov-norm} of Besov norm, we end up with
\begin{equation}\label{e:bd-G-u0}
\big\|\{(\partial_t \cg)_t -(\partial_t \cg)_s\} u_0\big\|_{\cb^{\ka,\mu_t}_{2,\infty}}\lesssim |t-s|^\ga \|u_0\|_{\ch^{1+\ka}}\, . 
\end{equation}
In the same way, we have
\begin{align*}
&\big\|\vp_{j}\ast (\partial_t \cg)_s u_0\big\|_{L^2_{\mu_s}} \\%\lesssim \big\|\big(\vp_j \ast((\partial_t G)_s \ast  u_0\big)\big\|_{L^2}\\
&\lesssim \bigg( \int_{\R^d}d\xi \, \big| \widehat{\vp_j}(\xi)  \big|^2\big|\cos(s|\xi|)\big|^2  \big|\widehat{u_0}(\xi)\big|^2\bigg)^{\frac12}\lesssim 2^{-j\ka} \bigg( \int_{\R^d}d\xi \, |\xi|^{2\ka} \big|\widehat{u_0}(\xi)\big|^2\bigg)^{\frac12}\lesssim 2^{-j\ka} \|u_0\|_{\ch^{1+\ka}}.
\end{align*}
Therefore, resorting to \eqref{e:Besov-norm} again, we get  $\sup_{s\in [0,T]} \big\|(\partial_t \cg)_s u_0\big\|_{\cb^{\ka,\mu_s}_{2,\infty}}\lesssim \|u_0\|_{\ch^{1+\ka}}$. Gathering this inequality with \eqref{e:bd-G-u0} we have thus shown that
\begin{equation}\label{e:bd-G-u0'}
\big\| (\partial_t \cg) u_0 \big\|_{\ce^{\ga,\ka}_{2,\infty}(T)} \lesssim \big\|u_0\big\|_{\ch^{1+\ka}}.
\end{equation}

Let us now handle the term $\cg_t u_1$ in \eqref{e:u-0}. It is treated similarly to $(\partial_t \cg)_t u_0$ in \eqref{e:bd-vp-G-u0}-\eqref{e:bd-G-u0'}, and we omit the details for sake of conciseness. Let us just mention that \eqref{e:bd-vp-G-u0} is replaced by  
\begin{align*}
\big\|\vp_{j}\ast \{\cg_t -\cg_s\} u_1\big\|_{L^2_{\mu_t}} 
&\lesssim \bigg( \int_{\R^d}d\xi \, \big| \widehat{\vp_j}(\xi)  \big|^2\frac{|\sin(t|\xi|)-\sin(s|\xi|)|^2}{|\xi|^2}  \big|\widehat{u_1}(\xi)\big|^2\bigg)^{\frac12}\\
&\lesssim 2^{-j\ka}|t-s| \bigg( \int_{\R^d}d\xi \, |\xi|^{2\ka} \big|\widehat{u_1}(\xi)\big|^2\bigg)^{\frac12}.
\end{align*}
Following the same steps as before, we then get 
\begin{equation}\label{e:bd-G-u1}
\big\| \cg u_1 \big\|_{\ce^{\ga,\ka}_{2,\infty}(T)} \lesssim \big\|u_1\big\|_{\ch^{\ka}}\,. 
\end{equation}
Plugging together \eqref{e:bd-G-u0'} and \eqref{e:bd-G-u1}, we have thus proved our claim \eqref{e:bd-u-0}. 
\end{proof}

With Lemma \ref{lem:ini-condi}  in hand, we now turn to the proof of our existence-uniqueness result.

\begin{proof}[Proof of Theorem \ref{theo:wellposedness}] Due to the expression \eqref{e:bd-u-0} for the initial condition $u^{(0)}$, we will not use a standard argument based on patching solutions defined on small intervals. Instead of that, our method will be based on Picard iterations.  In other words, we consider the sequence of processes $(u^{(\ell)})_{\ell \geq 0}$ defined as: for every $t\in [0,T]$, $u^{(0)}_t:=(\partial_t \cg)_t u_0+\cg_t u_1$ and
\begin{equation}\label{picard}
u^{(\ell+1)}_t:=(\partial_t \cg)_t u_0+\cg_t u_1+\int_0^t \cg_{t-r}(u^{(\ell)}_r \, d\dot W_r).
\end{equation}
By combining the results of Proposition \ref{prop:Riemann-sum} and Lemma \ref{lem:ini-condi}, we can immediately guarantee that $u^{(\ell)}$ is well defined in the space $\ce^{\ga,\ka}_{2,\infty}(T)$, for every $\ell \geq 0$. We now divide our proof into the existence and the uniqueness parts.

\noindent
\textit{Existence.} In order to show the convergence of $u^{(\ell)}$ in $\ce_{2,\infty}^{\gamma,\mu}(T)$, let us consider the difference $v^{(\ell)}:=u^{(\ell+1)}-u^{(\ell)}$.  Owing to \eqref{picard}, it is clear that for every $\ell \ge 0$ we have $v_0^{(\ell)}=0$ and $v^{(\ell)}$ satisfies 
\begin{equation}
v^{(\ell)}_t:=\int_0^t \cg_{t-r}(v^{(\ell-1)}_r \, d\dot W_r).
\end{equation}
We are here in a position to apply the estimate \eqref{boun-int-you-wav'} and assert that, for all $t\in [0,T]$
$$
\|v^{(\ell)}\|^q_{\ce^{\ga,\ka}_{2,\infty}(t)} \leq    c_T    \int_0^t \|v^{(\ell-1)}\|_{\ce^{\ga,\ka}_{2,\infty}(t_1)}^q \, dt_1\, ,
$$
where $c_T$ only depends on    $T$ and    $\|\dot W\|_{\cac^\theta([0,T];\cb^{-\al,P}_{p,\infty})}$, and $q$ is the (finite) parameter provided by Proposition \ref{prop:Riemann-sum}. Iterating this inequality, we deduce
$$
\|v^{(\ell)}\|_{\ce^{\ga,\ka}_{2,\infty}(T)}^q \leq (c_T)^\ell \int_0^T\int_0^{t_1}\ldots\int_0^{t_{\ell}} \|v^{(0)}\|_{\ce^{\ga,\ka}_{2,\infty}(t_{\ell})}^q \, dt_1\cdots dt_{\ell} \leq  \|v^{(0)}\|_{\ce^{\ga,\ka}_{2,\infty}(T)}^q \frac{(c_T T)^\ell}{\ell !},
$$
which can naturally be recast as
$$\|u^{(\ell+1)}-u^{(\ell)}\|_{\ce^{\ga,\ka}_{2,\infty}(T)} \leq \|u^{(1)}-u^{(0)}\|_{\ce^{\ga,\ka}_{2,\infty}(T)} \frac{\left(c_T T\right)^{\frac{\ell}{q}}}{(\ell !)^{\frac{1}{q}}}.$$
This proves that $(u^{(\ell)})$ is a Cauchy sequence in $\ce^{\ga,\ka}_{2,\infty}(T)$, and accordingly it converges to some limit $u$. Letting $\ell$ tend to infinity in \eqref{picard}, a standard procedure reveals  that $u$ is the desired solution to equation \eqref{equa-young-wellposed}.

\noindent
\textit{Uniqueness.} The argument is essentially the same as above. Namely, if $u,v$ are two solutions, we have for all $\ell\geq 1$,
$$
\|u-v\|_{\ce^{\ga,\ka}_{2,\infty}(T)}^q \leq (c_T)^\ell \int_0^T\int_0^{t_1}\ldots\int_0^{t_{\ell}} \|u-v\|_{\ce^{\ga,\ka}_{2,\infty}(t_{\ell})}^q \, dt_1\cdots dt_{\ell} \leq  \|u-v\|_{\ce^{\ga,\ka}_{2,\infty}(T)}^q \frac{(c_T T)^\ell}{\ell !},
$$
and we get the conclusion by letting $\ell$ tend to infinity.
\end{proof}

\section{Application to Gaussian noises}\label{sec:wave}
In this section, we prove that the Gaussian  noise whose covariance is defined by~\eqref{e:cov-sko} is an appropriate input for the wave equation. This amounts to prove that $\dot W$ sits in a space of the form $C^\theta([0,T], \cb_{p,\infty}^{-\alpha,P})$ as in Proposition \ref{prop:Riemann-sum}. Before proceeding to the proof, we first characterize the space $\cb_{p,\infty}^{-\alpha,P}$ thanks to simple properties of the Fourier transform. Let us start by defining a useful Fourier type operator. 

\begin{definition}\label{def:J}
Let $f$ be a smooth enough function defined on $\R^d$, and recall that its Fourier transform is denoted by $\cf f$. For an arbitrary constant $c>0$ and $s\in(-\infty, 1]$, we set
\[\cj^s_c f:=\cf^{-1}\big( \{1+c\,  |\cdot|^2\}^{\frac{s}{2}} \cf f\big)\, .\]
\end{definition}

We now upper bound norms in some weighted Besov spaces $\cb$ thanks to the Fourier operator $\cj$. This is summarized in the following Lemma.

\begin{lemma}\label{lem:new-norm} Consider $p\ge 2$ and $s\in(-\infty, 1]$. Recall that the spaces $L_w^p$ are given in Definition~\ref{def:L-mu} and the weighted Besov spaces $\cb_{p,q}^{s,\mu}$ are introduced in Definition \ref{def:Besov}. Let us also recall that we are working with the polynomial weight $P(x)=(1+|x|^{1+d})^{-1}$ and the related Besov spaces $\cb_{p,q}^{s,P}$ as in Lemma \ref{lem:change-weight}. Then there exists a constant $c_p>0$ depending on $p$    and $d$    only such that
\begin{equation}\label{new-norm}
\|f\|_{\cb^{s,P}_{p,\infty}} \lesssim \left\|\cj^s_{c_p} f \right\|_{L^p_P}\, ,
\end{equation}
where the operator $\cj_c^s$ has been introduced in Definition \ref{def:J}. 
\end{lemma}

\begin{proof} The proof of this Lemma is basically borrowed from \cite{rychkov}, to which we will refer for further details. First going back to the definition \eqref{e:Besov-norm} of Besov norms and invoking the fact that $\| \cdot \|_{\ell^\infty(\mathbb N)}\le \|\cdot\|_{\ell^p(\mathbb N)}$ for $1\le p<\infty$, we easily get that 
\[\|f\|_{\cb_{p,\infty}^{s,P}} \lesssim \|f\|_{\cb_{p,p}^{s,P}}\,.\]
Next, write equation \eqref{e:Besov-norm}, which gives 
\begin{align*}
&\|f\|_{\cb_{p,p}^{s, P}}^p = \sum_{j\geq 0} 2^{sj p} \int_{\R^{   d  }} dx \, P(x) |(\pmb{\vp}_j \ast f)(x)|^p\\
&=\int_{\R^{   d  }} dx \, P(x) \bigg(\sum_{j\geq 0} 2^{sj p}  |(\pmb{\vp}_j \ast f)(x)|^p\bigg)= \int_{\R^{   d  }} dx \, P(x) \left\|j\mapsto 2^{sj }  |(\pmb{\vp}_j \ast f)(x)| \right\|^p_{\ell^p(\mathbb N)}\,.
\end{align*}
Owing to the fact that $\|\cdot\|_{\ell^p(\mathbb N)}\le \|\cdot\|_{\ell^2(\mathbb N)}$ for $p\ge 2$, we obtain
\begin{equation}\label{e:f-B-norm}
\|f\|_{\cb_{p,p}^{s, P}}^p 
\leq \int_{\R^{   d  }} dx \, P(x) \bigg(\sum_{j\geq 0} 2^{2sj }  |(\pmb{\vp}_j \ast f)(x)|^2\bigg)^{\frac{p}{2}} \, .
\end{equation}
Note that the right-hand side of \eqref{e:f-B-norm} is exactly $\|f\|^p_{F_{p,2}^{s,P}}$ (see (2.2) in \cite{rychkov} for the definition of $\|\cdot\|_{F_{p,q}^{s,P}}$). Hence \eqref{e:f-B-norm} can be rephrased as
\[\|f\|_{\cb^{s,P}_{p,p}} \leq \|f\|_{F^{s,P}_{p,2}} \, .\]

The estimate \eqref{new-norm} now stems from the combination of \cite[Theorem 2.18]{rychkov} and \cite[Theorem 1.10]{rychkov}, which allows us to assert that for some suitable $c_p>0$,
$$\big\|f\big\|_{F^{s,P}_{p,2}} \lesssim \big\|\cj^s_{c_p} f \big\|_{L^p_P}\, ,$$
for some proportional constant that only depends on $p$    and $d$  .
\end{proof}

We now turn to a definition of our Gaussian noise allowing some proper couplings for approximations. Namely we will define $\dot{W}$ through an harmonizable representation of the form
\begin{equation}\label{e:dot-W}
\dot W_t(x):=c_{a_0}\int_{\la\in \R}\int_{\beta\in \R^d} \, \frac{e^{\imath \la t}-1}{|\la|^{\frac{3-a_0}2}}  \,  e^{\imath\beta\cdot x}   \,\widehat{B}(d\la,d\beta) \, .
\end{equation}
Here we use a parameter $a_0\in(0,2)$, $c_{a_0}$ is some proper positive constant,  and $\widehat B$ is a complex-valued Gaussian random measure  on $\mathcal B(\R\times \R^d)$ such that for any $I_1, I_2 \in \mathcal B(\R)$ and $A_1,A_2\in\mathcal B(\R^d)$, 
\begin{equation}\label{k1}
\E[\hat B(I_1\times A_1)]=0, \quad \E[\hat B(I_1\times A_1)\overline{\hat B(I_2\times A_2)}] 
=
m(I_1\cap J_1) \mu(A_1\cap A_2),
\end{equation}
where $m$ is the Lebesgue measure on $\R$ and $\mu(d\beta)$ is the spectral measure of a nonnegative definite (generalized) function $\gamma(x)$. In order to get a real-valued noise, we also assume that the Gaussian measure $\widehat{B}$ is such that
\begin{equation}\label{k2}
\overline{\hat B(I_1\times A_1)}=\hat B\lp -\lp I_1\times A_1\rp\rp \, .
\end{equation}
It is readily checked from \eqref{k1}-\eqref{k2} that for an appropriate value of $c_{a_0}$, the covariance of $\dot{W}$ is (formally) given by 
\begin{equation}\label{e:cov-general}
\E[\dot W_t(x) \dot W_s(y)]=  R_{a_0}(t,s)\gamma(x-y), 
\end{equation}
with
$$R_{a_0}(t,s):=\frac12\big\{|t|^{2-a_0}+|s|^{2-a_0}-|t-s|^{2-a_0}\big\} .$$

%\textbf{[Aur] A possible minor issue in the paper: we use the same notation $\dot W$ for the space-time derivative and the space (only) derivative of $W$. I don't know if it's worth fixing.}

The definition \eqref{e:dot-W} of $\dot{W}$ comes with a natural approximation by smooth functions. Namely for $n\ge 1$ we define
\begin{equation}\label{e:dot-W-n'}
\dot W^n_t(x):=c_{a_0}\int_{\la\in \R}\int_{\beta\in \mathbf B_{2^n}}\widehat{B}(d\la,d\beta) \, \frac{e^{\imath \la t}-1}{|\la|^{\frac{3-a_0}2}} \,  e^{\imath \beta\cdot x}   \, ,
\end{equation}
where we recall that the balls $\mathbf{B}_{R}$ are introduced in Notation \ref{notations}.
We shall show that $\dot W^n$ converges to $\dot W_t(x)$ given in \eqref{e:dot-W} in a proper space.

\begin{theorem}\label{prop:regu-frac-noise}
 Let    $a_0\in(0,2)$ and consider $\theta\in(0, 1-\frac{a_0}2)$. Let $\{\dot W^n, t\in[0,T], x\in\R\}$ be the field defined by \eqref{e:dot-W-n'}. We assume  the existence of a threshold parameter $\al_d>0$ such that the measure $\mu$ satisfies
\begin{equation}\label{e:con-gamma}
\int_{\R^d} \left(\frac{1}{1+|\beta|} \right)^{2\alpha_d+\varepsilon}\mu(d\beta)<\infty\,, ~\text{ for every } ~\varepsilon>0\,.
\end{equation}
 Then for all $\al>\al_d$  and $p$ large enough,  $\{\dot W^n,n \ge1\}$ forms a Cauchy sequence in the space
 $$
 L^p\lp \Omega; \cac^{\theta}\lp [0,T], \cb^{-\al,P}_{p,\infty}\rp\rp ,
 $$
  with the limit denoted by $\dot W$.    The noisy input $\dot W$ satisfies the following properties:
\begin{enumerate}[wide, labelwidth=!, labelindent=0pt, label=\textnormal{(\roman*)}]
\setlength\itemsep{.02in}

\item 
Almost surely we have $\dot W\in \cac^\theta([0,T], \cb_{p,\infty}^{-\alpha, P})$.
\item 
The covariance function of $\dot W$ is given by \eqref{e:cov-general}. 
\end{enumerate}
\end{theorem}

\black 

\begin{proof}
In the sequel we will write  $\delta \dot W_{st}^n = \dot W_t^n-\dot W_s^n$ to alleviate notations. We will reduce our proof to moment computations thanks to two standard steps:\\
(1) Invoking a telescoping sum argument, the convergence of the sequence $\{\dot W^n; n\ge 1\}$ is implied by the convergence of the series 
\begin{equation}\label{e:series-1'}
\sum_{n=1}^\infty \left(\E\left[\big\|\dot W^{n+1}-\dot W^n\big\|^p_{\cac^\theta([0,T], \cb^{-\al,P}_{p,\infty})} \right]\right)^{1/p}\,.
\end{equation}
(2) Owing to a standard application of Garsia's lemma the following holds true: if we assume that for each summand in \eqref{e:series-1'} and for $p\ge \frac{1}{1-\frac{a_0}2-\theta}$, there exists $\varepsilon>0$ such that
\begin{equation}\label{e:bd-delta-W'}
\E\left[\left\|\delta \dot W_{st}^{n+1}  - \delta \dot W_{st}^{n} \right\|^p_{\cb_{p,\infty}^{-\alpha, P}} \right]\lesssim \frac{|t-s|^{(1-\frac{a_0}{2}) p}}{2^{\varepsilon n p }}\,, 
\end{equation}
then we also have 
\begin{equation}\label{e:bd-W'}
\E\left[\big\|\dot W^{n+1}-\dot W^n\big\|^p_{\cac^\theta([0,T], \cb^{-\al,P}_{p,\infty})} \right]\lesssim \frac{1}{2^{\varepsilon n p}}\, .
\end{equation}
Putting together \eqref{e:series-1'} and \eqref{e:bd-W'}, we are now reduced to prove \eqref{e:bd-delta-W'}.

In order to upper bound the right-hand side of \eqref{e:bd-delta-W'}, we invoke \eqref{new-norm} which can be read  
\begin{multline}
\E\left[\left\|\delta \dot W_{st}^{n+1}  - \delta \dot W_{st}^{n} \right\|^p_{\cb_{p,\infty}^{-\alpha, P}} \right]  \\
\lesssim 
\int_{  \R^{d+1}  } \frac{dx}{1+|x|^{   1+d  }} \, \mathbb{E}\Big[\Big| \cf^{-1}\Big(\{1+|\cdot |^2\}^{-\frac{\al}{2}}\cf\big(\delta\dot{W}_{st}^{n+1}-\delta\dot{W}_{st}^n\big)\Big)(x)\Big|^{p}\Big]\,.
\end{multline}
In addition, $\dot W^{n+1}-\dot W^n$ is a Gaussian process and $\cf$ is a linear transform. Hence we get 
\begin{equation}\label{e:bd-delta-W-1'}
\E\left[\left\|\delta \dot W_{st}^{n+1}  - \delta \dot W_{st}^{n} \right\|^p_{\cb_{p,\infty}^{-\alpha, P}} \right]\lesssim \int_{\R^{   d   }} \frac{dx}{1+|x|^{   1+d  }} \,\big( Q_{st}^n (x) \big)^{p/2}\,,
\end{equation}
where we have set 
\begin{equation*}
Q_{st}^{n}(x)
:=
\mathbb{E}\Big[\Big| \cf^{-1}\Big(\{1+|\cdot|^2\}^{-\frac{\al}{2}}\cf\big(\delta\dot W^{n+1}_{st}-\delta \dot W^n_{st}\big)\Big)(x)\Big|^2\Big]\, .
\end{equation*}
Let us focus our attention on the evaluation of $Q^n$ above, defined for $0\le s<t\le T$ and $x\in\R^{   d  }$. To this aim, by simply writing the definition of Fourier transform we get
\begin{multline}\label{e:Q-st'}
Q^n_{st}(x) = \int_{\R^{   d  }}\int_{\R^{   d  }} d\xi d\xiti \, e^{\imath x\cdot (\xi-\xiti)} \{1+|\xi|^2\}^{-\frac{\al}{2}} \{1+|\xiti|^2\}^{-\frac{\al}{2}} \\ 
\times \int_{\R^{   d  }}\int_{\R^{   d  }} dy d\yti \, e^{-\imath \xi\cdot  y}e^{\imath \xiti\cdot \yti}  
\, \mathbb{E}\Big[\big\{\delta\dot{W}^{n+1}_{st}-\delta\dot{W}^n_{st}\big\}(y) 
\overline{\big\{\delta \dot{W}^{n+1}_{st}-\delta\dot{W}^n_{st}\big\}(\yti)}\Big]\,.
\end{multline}
In addition, invoking relation \eqref{e:dot-W-n'} and the covariance structure of a complex valued white noise, we obtain
\begin{equation*}
\mathbb{E}\Big[\big\{\delta\dot{W}^{n+1}_{st}-\delta\dot{W}^n_{st}\big\}(y) 
\overline{\big\{\delta \dot{W}^{n+1}_{st}-\delta\dot{W}^n_{st}\big\}(\yti)}\Big]=c\, |t-s|^{2-a_0} \int_{   \bfb_{2^{n+1}}\backslash \bfb_{2^n}  }  e^{\imath \beta\cdot( y-\yti)}  \mu(d\beta)      \, , 
\end{equation*}  
for some constant $c$. Plugging this inequality into \eqref{e:Q-st'} and carefully computing the Fourier transforms reveals that \begin{align}\label{e:Q'}
Q^n_{st}(x)&=c \, |t-s|^{2-a_0} \int_{   \bfb_{2^{n+1}}\backslash \bfb_{2^n}  } \frac{ \mu(d\beta)  }{\{1+|\beta|^2\}^{\al}}   \,.
\end{align}

We now bound the right-hand side of \eqref{e:Q'}:  by setting $\varepsilon:=\al-\al_d>0$, we get
\begin{align}
\int_{   \bfb_{2^{n+1}}\backslash \bfb_{2^n}  } \frac{\mu(d\beta)}{\{1+|\beta|^2\}^{\al}} & = \int_{   \bfb_{2^{n+1}}\backslash \bfb_{2^n}  } \{1+|\beta|^2\}^{-\frac\varepsilon2}\frac{\mu(d\beta)}{\{1+|\beta|^2\}^{\al_d+\frac{\varepsilon}{2}}}  
\lesssim 2^{-n\varepsilon}   ,\label{e:bd-int'}
\end{align}
where the inequality follows from the assumption \eqref{e:con-gamma}. 
Reporting \eqref{e:bd-int'} into \eqref{e:Q'}, we have thus obtained 
\[
Q_{st}^n(x) \lesssim 2^{-n\varepsilon} |t-s|^{2-a_0}\,.
\]
   Plugging this information into \eqref{e:bd-delta-W-1'}, this trivially yields our claim \eqref{e:bd-delta-W'}. Since we have argued that the convergence of $\dot W^n$ could be reduced to \eqref{e:bd-delta-W'}, we have proved that $\dot W^n$ converges in $L^p(\Omega; \cac^{\theta}([0,T], \cb^{-\al,P}_{p,\infty}))$ for $p\ge \frac1{1-\frac{a_0}{2}-\theta}$. 
In particular, this proves item $(i)$ in our Theorem. 

Finally, item $(ii)$ in our Theorem can be proved by a standard procedure on covariance functions from \eqref{e:dot-W-n'} to \eqref{e:cov-general}. 
The proof is concluded. 
\end{proof}

\begin{remark}\label{remark:special-cases}
For the special cases $\gamma(x)=\prod_{i=1}^d |x_i|^{-a_i}$ and $ \gamma(x)=|x|^{-a}$ with $a_i\in(0,2), a\in(0,2d)$,\footnote{Here we use the convention that for $x\in \R$,  $|x|^{-1}:=\delta(x)$ where $\delta(x)$ is the Dirac delta function and $|x|^{-a}:=(|x|^{-a+2})''$ for $a\in(1,2)$ where the second derivative is taken in the distributional sense. Similar convention also applies for $x\in\R^d$.} the spectral measures are  $\mu(d\beta)=\prod_{i=1}^d |\beta_i|^{-(1-a_i)} \, d\beta$ and $\mu(d\beta)=|\beta|^{-(d-a)} \, d\beta$, respectively. If we use the notation $a=\sum_{i=1}^d a_i$, condition \eqref{e:con-gamma} is verified for the optimal threshold
$\al_d = \frac{a}{2}$. 
\end{remark}

\begin{remark}
In the current paper, we have chosen  the harmonic representation~\eqref{e:dot-W-n'} to approximate the Gaussian noise, which is convenient in our setting. We could also have used, similarly to what is done in \cite{cdot}, a sequence  of mollifications $\dot W^n=\vp_n*\dot W$. This second approach might be closer in spirit to our previous work \cite{CDST}, but is slightly more cumbersome.  
\end{remark}

Let us recall that we have denoted by $\rho_d$ the regularization parameter exhibited in Proposition \ref{prop:regu-effect}, with explicit value
$$
\rho_d:=\begin{cases}
1 & \text{if} \ d=1\\
\frac12 & \text{if} \ d=2.
\end{cases}
$$

\begin{proposition} \label{prop:wave-eq}
   Assume that $d\in \{1,2\}$.    Under the same setting as for Theorem \ref{prop:regu-frac-noise}, assume that  $a_0<\rho_d$ and that condition \eqref{e:con-gamma} is satisfied for $\al_d$ such that
   \begin{equation}\label{e:con-sol}
0<\al_d <\frac12 (\rho_d-a_0).
\end{equation}
 Then there exist coefficients $0<\ka,\al,\ga,\theta <1$ and $p>   d+1  $ that satisfy~\eqref{cond-ka-al-ga-the} and such that $\dot{W}\in \cac^\theta([0,T];\cb^{-\al,P}_{p,\infty})$ almost surely. In particular Theorem \ref{theo:wellposedness} can be applied to the noise $\dot W$ and equation \eqref{equa-young-wellposed} admits a unique solution in $\ce_{2,\infty}^{\gamma, \kappa}(T)$. 
\end{proposition}
\begin{proof} 
Thanks to condition \eqref{e:con-sol}, we can fix $\varepsilon >0$ small enough so that $\al_d <\frac12 (\rho_d-a_0)-2\varepsilon$, $\theta:=1-\frac{a_0}2-\varepsilon   \in [0,1]$ and $\al:=\al_d+\varepsilon\in [0,1]$. By Theorem \ref{prop:regu-frac-noise}, we know that $\dot W\in C^\theta([0,T]; \cb_{p,\infty}^{-\alpha, P})$. Besides, observe that
\[ \theta-\alpha=1-\frac{a_0}2-\al_d-2\varepsilon > 1-\frac{\rho_d}{2}\,.\]
We can thus pick two parameters $p$ large enough and $\delta>0$ small enough so that
\begin{equation}\label{e:theta-alpha}
\theta-\alpha>   1-\frac{\rho_d}{2}  +\frac{   d  }{2p}+\delta\,.
\end{equation}

 We now consider $\kappa, \gamma$ such that the second and third conditions in \eqref{cond-ka-al-ga-the} are fulfilled, of the form 
 \begin{equation}\label{e:ka-ga}
 \kappa=\alpha+\frac{   d  }{p}+\delta, \mbox{ and } \gamma=1-\theta+\delta\,,
 \end{equation}
for $\delta>0$ small enough and $p\ge 1$ large enough. Notice that since we have $\alpha, \theta\in[0,1]$, we also have $\kappa, \gamma\in[0,1]$ whenever $p$ is large and $\delta$ is small. Moreover, for $\gamma$ as in \eqref{e:ka-ga}, the last condition in \eqref{cond-ka-al-ga-the} is verified:
\[1-\frac{   d+1  }{p}>\gamma\,.\]

It remains to prove that the first condition in \eqref{cond-ka-al-ga-the} holds true. Towards this aim, consider $\kappa, \gamma$ as in \eqref{e:ka-ga} and compute 
\begin{equation}\label{e:ka-al-ga}
  \kappa+\alpha+\gamma+(1-\rho_d)=2\alpha+\frac{d}{p}+1-\theta+2\delta+(1-\rho_d)= 2-\rho_d-2(\theta-\alpha)+\frac{d}{p}+2\delta+\theta\,.  
\end{equation}
In addition, resorting to \eqref{e:theta-alpha}, we have $2(\theta-\alpha)>   2-\rho_d+\frac{d}{p}+2\delta $. Plugging this relation into~\eqref{e:ka-al-ga}, we get 
\[\kappa+\alpha+\gamma   +(1-\rho_d)  <\theta\,.\] 
This shows that the set of conditions \eqref{cond-ka-al-ga-the} is satisfied and that $\dot W$ is an element of $C^{\theta}([0,T]; \cb_{p,\infty}^{-\alpha, P})$. Our proof is achieved. 
\end{proof}

Our    final    result draws on Proposition \ref{prop:wave-eq} and gives conditions on    $a_0, a$ when the covariance function $\gamma(x)$ satisfies the scaling property \eqref{e:gamma-scaling}.

\begin{corollary} \label{cor:wave-eq}
   Assume that $d\in \{1,2\}$.    Under the same setting as in Remark \ref{remark:special-cases}, assume that    
\begin{equation}\label{e:con-special}
 a_0+a<\rho_d\,.
   \end{equation}
    Then there exist coefficients $0<\ka,\al,\ga,\theta <1$ and $p>   d+1  $ that satisfy~\eqref{cond-ka-al-ga-the} and such that $\dot{W}\in \cac^\theta([0,T];\cb^{-\al,P}_{p,\infty})$ almost surely. In particular Theorem \ref{theo:wellposedness} can be applied to the noise $\dot W$ and equation \eqref{equa-young-wellposed} admits a unique solution in $\ce_{2,\infty}^{\gamma, \kappa}(T)$. 
\end{corollary}
\begin{proof}
 As noticed in Remark \eqref{remark:special-cases}, condition \eqref{e:con-gamma} is satisfied by $\al_d=\frac{a}{2}$ in this situation. The desired result  then  follows from a  straightforward verification that \eqref{e:con-special} is equivalent to \eqref{e:con-sol}  when $\al_d=\frac{a}{2}$.
\end{proof}

\appendix

\section{}\label{sec:appendix}

This Appendix summarizes some technical computations which play a prominent role in our Besov space considerations. 

\subsection{A local reproducing formula}

We start by labelling a Theorem borrowed from \cite[Theorem 1.6]{rychkov}:

\begin{theorem}[Rychkov]\label{theo:local-reproducing-formula}
Let $\vp_0\in \cac_c^\infty$ be such that $\int_{\R^d}\vp_0(x)\, dx \neq 0$. Then for any integer $L\geq 0$, there exist two functions $\psi_0,\psi\in \cac_c^\infty$ such that:

\begin{enumerate}[wide, labelwidth=!, labelindent=0pt, label=(\roman*)]
\setlength\itemsep{.02in}

\item
The function $\psi$ is an element of  $\cd_L$,  where we recall that $\cd_L$ is introduced in Definition~\ref{def:D-L}

\item
For every distribution $f\in \cd'$, the following decomposition holds true in $\cd'$:
\begin{equation*}
f=\sum_{j\geq 0} \vp_j \ast \psi_j \ast f \, .
\end{equation*}
\end{enumerate}
\end{theorem}

\subsection{Proof of Lemma \ref{lem:multipl}} \label{sec:proof-lem}

Recall that $\pmb{\vp}_0,\pmb{\vp}_j$ stand for the test-functions introduced in Notation \ref{notation:pmb-vp}, and consider an arbitrary number $L\geq 2$. By Theorem \ref{theo:local-reproducing-formula}, there exist two functions $\psi_0,\psi\in \cac_c^\infty$ such that $\psi\in \cd_L$ and for every $h\in \cd'$,
\begin{equation}\label{repro-prod}
h=\sum_{\ell\geq 0} \psi_\ell\ast \pmb{\vp}_\ell\ast h \, .
\end{equation}
\textit{As in Lemma \ref{lem:vp-j-ast-psi-l}, for the sake of clarity, we will assume in the sequel that the supports of $\pmb{\vp}_0$, $\pmb{\vp}$, $\psi$ and $\psi_0$ are all included in the interval $\mathbf B_1$, but the arguments could be easily extended to any supporting interval $\mathbf B_K$}.

For all $\ell,k\geq 0$, let us set $\delta_\ell f:=\psi_\ell \ast \pmb{\vp}_\ell\ast f$ and $\delta_k g:=\psi_k \ast \pmb{\vp}_k\ast g$. Owing to formula \eqref{repro-prod}, one can proceed similarly to \eqref{defi-i-ii} and split $\vp_j*(f\cdot g)$ as 
\begin{align}\label{e:vp*fg}
\pmb{\vp}_j \ast (f\cdot g)&=\sum_{\ell,k\geq 0} \pmb{\vp}_j \ast (\delta_\ell f\cdot \delta_k g) =: \ I_j+I\! I_j\, ,
\end{align}
where the terms $I_j$ and $I\! I_j$ are defined by
\begin{equation}\label{e:I-II}
I_j =\sum_{\substack{\ell,k\geq 0\\0\leq \ell\leq \max(j,k) }} \pmb{\vp}_j \ast (\delta_\ell f\cdot \delta_k g)\,, \quad I\!I_j=\sum_{\substack{\ell,k\geq 0\\ \ell> \max(j,k) }} \pmb{\vp}_j \ast (\delta_\ell f\cdot \delta_k g)\,.  
\end{equation}
We now bound $I_j$ and $I\!I_j$ in two different steps. 

\smallskip

\noindent
\textit{Treatment of $I_j$}. We first bound each term $\pmb{\varphi}_j*(\delta_\ell f\cdot \delta_kg)$ in the sum defining $I_j$ in \eqref{e:I-II}. That is, resorting to \eqref{e:conv-Young} together with the fact that $\pmb{\vp}_j$  is bounded in $L^1$, we get 
\[\|\pmb{\vp}_j*(\delta_\ell f\cdot \delta_k g) \|_{L_\mu^p}\lesssim \|\delta_\ell f\cdot \delta_k g \|_{L_\mu^p}\,.\]
Let us now consider $p_1, p_2, \mu_1, \mu_2$ as in \eqref{e:con-parameters}. A standard application of H\"older's inequality shows that 
\[ \|\delta_\ell f\cdot \delta_k g \|_{L_\mu^p}\le \|\delta_\ell f\|_{L_{\mu_1}^{p_1}} \| \delta_k g \|_{L_{\mu_2}^{p_2}}\,.\]
Next for each block $\delta_\ell f,  \delta_k g$,  we use the fact that $\psi_\ell$ is bounded in $L^1$ (uniformly in $\ell$), the expression \eqref{e:Besov-norm} for the norm in Besov spaces,  and Young's inequality. We get 
\begin{equation}\label{e:del-f-g}
\|\delta_\ell  f \cdot \delta_k g\|_{L_\mu^p} \lesssim 2^{\alpha\ell}2^{-\beta k} \|f\|_{\cb^{-\al,\mu_1}_{p_1,\infty}}\|g\|_{\cb^{\beta,\mu_2}_{p_2,\infty}}\,.
\end{equation}
Summarizing our considerations so far, we have obtained 
\[\|\pmb{\vp}_j*(\delta_\ell f\cdot \delta_k g) \|_{L_\mu^p}\lesssim  2^{\alpha\ell}2^{-\beta k} \|f\|_{\cb^{-\al,\mu_1}_{p_1,\infty}}\|g\|_{\cb^{\beta,\mu_2}_{p_2,\infty}}\,, \]
for a proportional constant that only depends on $  \mu_\ast  $. Reporting this estimate into the sum \eqref{e:I-II} defining $I_j$, we end up with 
\begin{equation}\label{e:I-j}
\big\|I_j\big\|_{L^p_\mu}  \lesssim \|f\|_{\cb^{-\al,\mu_1}_{p_1,\infty}}\|g\|_{\cb^{\beta,\mu_2}_{p_2,\infty}} \sum_{\substack{\ell,k\geq 0\\0\leq \ell\leq \max(j,k) }}  2^{\al\ell} 2^{-\beta k}\,.
\end{equation}
Let us upper bound the sum in the right-hand side of \eqref{e:I-j}. Splitting the summation in a straightforward way, we get 
\begin{align*}
\sum_{\substack{\ell,k\geq 0\\0\leq \ell\leq \max(j,k) }}  2^{\al\ell} 2^{-\beta k}
&\lesssim  \sum_{0\leq k\leq j}2^{-\beta k}\sum_{0\leq \ell\leq j}  2^{\al\ell} +\sum_{k>j}2^{-\beta k}\sum_{0\leq \ell\leq k}  2^{\al\ell} \\
&\lesssim 2^{\al j} \sum_{k\geq 0}2^{-\beta k} +\sum_{k\geq 0}2^{-(\beta-\al) k}
\lesssim 2^{\al j}+1 \lesssim 2^{\al j}\, ,
\end{align*}
where we have used the fact that $\beta>\alpha$ (see \eqref{e:con-parameters}). Going back to \eqref{e:I-j}, we have thus obtained 
\begin{equation}\label{e:I-j'}
\|I_j\|_{L_\mu^p}\lesssim 2^{\alpha j} \|f\|_{\cb^{-\al,\mu_1}_{p_1,\infty}}\|g\|_{\cb^{\beta,\mu_2}_{p_2,\infty}}\,.
\end{equation}

\noindent
\textit{Treatment of $I\! I_j$}. As in the previous step, we first analyze a single term $\pmb{\vp}_j*(\delta_\ell f\cdot \delta_k g)$. In the regime $\ell >\max(j, k)$, we write the definition of convolution product and apply Fubini in order to get
\begin{align}
\big(\pmb{\vp}_j \ast \big[ \delta_\ell f \cdot \delta_k g\big]\big)(x)&=\int dy \, \pmb{\vp}_j(x-y) \big[ \psi_\ell \ast (\pmb{\vp}_\ell\ast f)\big](y) \big[ \psi_k \ast (\pmb{\vp}_k \ast g)\big](y)\nonumber\\
&=\int dz \int dv \, (\pmb{\vp}_\ell \ast f)(z) (\pmb{\vp}_k \ast g)(v) K_{jk\ell}(x, v, z)\, ,\label{term-ii-j}
\end{align}
where we have set 
\begin{equation}\label{e:K-j-k-l}
K_{jk\ell}(x, v, z):=\int dy \, \pmb{\vp}_j(x-y)\psi_\ell(y-z) \psi_k(y-v)\,.
\end{equation}

Let us further investigate the term $K_{jk\ell} $ defined in \eqref{e:K-j-k-l}. Resorting to a change of variable $y:=y-z$ and expressing the scalings in the functions $\pmb{\vp}$ and $\psi$, we have 
 \begin{align*}
K_{jk\ell}(x, v, z)&=\int dy \, \psi_\ell(y) \pmb{\vp}_j((x-z)-y)\psi_k(y+(z-v))\\
&=2^{d(\ell+j+k)}\int dy \, \psi(2^\ell y) \pmb{\vp}(2^j(x-z)-2^j y)\psi(2^k y+2^k(z-v))\\
&=2^{d(j+k)}\int dy \, \psi(y) \pmb{\vp}(2^j(x-z)-2^{-(\ell-j)} y)\psi(2^{-(\ell-k)} y+2^k(z-v))\,,
\end{align*}
where the last identity stems from another elementary change of variable $y:=2^{\ell} y$. For $m\ge 0$, we now introduce the indicator function $\chi_m:=\1_{\mathbf B_{2^{-m+1}}}$. Thanks to the fact that the support of both $\pmb{\vp}$ and $\psi$ is a subset of $\mathbf B_1$, for $\ell >\max(j,k)$, we end up with 
\begin{multline}\label{e:K-j-k-l-1}
K_{jk\ell} (x,v,z) = 2^{d(j+k)} \chi_j(x-z)\chi_k(z-v) \\
\times \int dy \, \psi(y) \pmb{\vp}(2^j(x-z)-2^{-(\ell-j)} y)\psi(2^{-(\ell-k)} y+2^k(z-v))\,.
\end{multline}
We now proceed as for \eqref{e:phi*psi}-\eqref{e:est1}, invoking the fact that $\psi\in \cd_L$ and using a Taylor expansion for $y\mapsto \pmb{\vp}(2^j(x-z)-2^{-(\ell-j)} y)\psi(2^k(z-v)+2^{-(\ell-k)} y)$. We let the reader check that we get 
$$\bigg|\int dy \, \psi(y) \pmb{\vp}(2^j(x-z)-2^{-(\ell-j)} y)\psi(2^k(z-v)+2^{-(\ell-k)} y)\bigg|\lesssim 2^{-L(\ell-\max(j,k))}\, ,$$
for some proportional constant that does not depend on $j,k,\ell$. Taking the above inequality into account in \eqref{e:K-j-k-l-1}, this leads to 
\begin{equation}\label{e:K-j-k-l-2}
|K_{jk\ell} (x,v,z)|\lesssim 2^{d(j+k)-L(\ell-\max(j,k)) }\chi_j(x-z)\chi_k(z-v)\,.
\end{equation}

Let us go back to \eqref{term-ii-j} with \eqref{e:K-j-k-l-2} in hand. This allows to write 
\begin{equation}\label{term-ii-j-1}
\Big|\big(\pmb{\vp}_j \ast \big[ \delta_\ell f \cdot \delta_k g\big]\big)(x)\Big|\lesssim 2^{-L(\ell-\max(j,k))+d(j+k)} \tilde K_{jk\ell}(x)\,,
\end{equation}
where the function $\tilde K_{jk\ell}$ can be expressed as 
\[\tilde K_{jk\ell}(x)= \int dz \int dv \, |(\pmb{\vp}_\ell \ast f)(z)| |(\pmb{\vp}_k \ast g)(v)| \chi_j(x-z)\chi_k(z-v)\,.\]
We now express $\tilde K_{jk\ell}$ as a convolution product thanks to elementary changes of variable. Specifically, setting $z:=x-z$ and $v=x-z-v$, we get 
\[\tilde K_{jk\ell} (x) = \int dz\, |(\pmb{\vp}_\ell \ast f)(x-z)| \chi_j(z)\int dv \,  |(\pmb{\vp}_k \ast g)(x-z-v)|  \chi_k(v)\,.\]
Hence it is easily seen that 
\begin{align*}
\tilde K_{jk\ell}(x) &=  \int dz\, |(\pmb{\vp}_\ell \ast f)(x-z)| \chi_j(z)\big[\chi_k\ast |\pmb{\vp}_k\ast g| \big](x-z)\\
&=  \bigg[\chi_j \ast \Big[|\pmb{\vp}_\ell \ast f| \cdot \big[ \chi_k\ast |\pmb{\vp}_k\ast g|\big]\Big] \bigg](x)\, .
\end{align*}

We plug the above formula into \eqref{term-ii-j-1}. Then our expression with convolutions allows the application of \eqref{e:conv-Young}, as in the previous step. We first get 
\[\big\|\pmb{\vp}_j \ast \big[ \delta_\ell f \cdot \delta_k g\big]\big\|_{L^p_\mu}\lesssim 2^{-L(\ell-\max(j,k))}\big(2^{dj} \|\chi_j\|_{L^1}\big)\, 2^{dk} \Big\| |\pmb{\vp}_\ell \ast f| \cdot \big[ \chi_k\ast |\pmb{\vp}_k\ast g|\big]\Big\|_{L^p_\mu}\,.
\]
Then thanks to H\"older's inequality, we obtain
\[\big\|\pmb{\vp}_j \ast \big[ \delta_\ell f \cdot \delta_k g\big]\big\|_{L^p_\mu}\lesssim 2^{-L(\ell-\max(j,k))}\big\|\pmb{\vp}_\ell \ast f\big\|_{L^{p_1}_{\mu_1}}\big(2^{dk}  \big\| \chi_k\ast |\pmb{\vp}_k\ast g|\big\|_{L^{p_2}_{\mu_2}}\big) \,.\]
Then owing to \eqref{e:conv-Young} again and similarly to \eqref{e:del-f-g}, this yields 
\begin{align}\label{e:vp-f-g}
\big\|\pmb{\vp}_j \ast \big[ \delta_\ell f \cdot \delta_k g\big]\big\|_{L^p_\mu}
\lesssim 2^{-L(\ell-\max(j,k))} 2^{\al \ell} 2^{-\beta k} \|f\|_{\cb^{-\al,\mu_1}_{p_1,\infty}}\|g\|_{\cb^{\beta,\mu_2}_{p_2,\infty}}\, ,
\end{align}
where the proportional constants only depend on $  \mu_\ast  $.

We can now proceed as in \eqref{e:I-j}-\eqref{e:I-j'}. That is we plug \eqref{e:vp-f-g}    into    the definition \eqref{e:I-II} of $I\!I_j$. This yields 
\begin{equation}\label{e:II-j}
\|I\!I_j\|_{L_\mu^p} \lesssim \|f\|_{\cb^{-\al,\mu_1}_{p_1,\infty}}\|g\|_{\cb^{\beta,\mu_2}_{p_2,\infty}}  \sum_{\substack{\ell,k\geq 0\\\ell > \max(j,k)}} 2^{-L(\ell-\max(j,k))}2^{\alpha\ell}2^{-\beta k}\,. 
\end{equation}
Then we split the sum in the right-hand side above and take into account the fact that $
\beta>\alpha, L>\alpha$ in order to obtain 
\begin{align*}
 \sum_{\substack{\ell,k\geq 0\\\ell > \max(j,k)}} 2^{-L(\ell-\max(j,k))}2^{\alpha\ell}2^{-\beta k}&= \sum_{0\leq k\leq j}2^{-\beta k} \sum_{\ell >j} 2^{-L(\ell-j)} 2^{\al \ell} +\sum_{k>j}2^{-\beta k}\sum_{\ell >k}2^{-L(\ell-k)} 2^{\al \ell}\\
&= 2^{\al j}\sum_{   0\leq k\leq j  }2^{-\beta k} \sum_{\ell >0} 2^{-(L-\al)\ell}  +\sum_{   k>j  }2^{-(\beta-\al) k}\sum_{\ell >0}2^{-(L-\al)\ell}  \\
&\lesssim    2^{\al j}+2^{-(\beta-\al)j}\lesssim 2^{\al j}\, .  
\end{align*}
Resorting to this inequality in \eqref{e:II-j}, we have obtained 
\begin{equation}\label{e:II-j'}
\|I\!I_j\|_{L_\mu^p} \lesssim  2^{\alpha j} \|f\|_{\cb^{-\al,\mu_1}_{p_1,\infty}}\|g\|_{\cb^{\beta,\mu_2}_{p_2,\infty}} \,. 
\end{equation}

\noindent\textit{Conclusion.} We simply gather \eqref{e:I-j'} and \eqref{e:II-j'} into \eqref{e:I-II} and \eqref{e:vp*fg}. This yields
\[\|\vp_j*(f\cdot g) \|_{L_{\mu}^p} \lesssim 2^{\alpha j} \|f\|_{\cb^{-\al,\mu_1}_{p_1,\infty}}\|g\|_{\cb^{\beta,\mu_2}_{p_2,\infty}} \,. \]
Owing to the definition \eqref{e:Besov-norm} of norms in Besov spaces, we thus easily get the fact that $   f\cdot g  $ sits in $\cb_{p,\infty}^{-\alpha, \mu}$ and that  \eqref{e:product} holds. This finishes the proof.  \qed

\subsection{Proof of Lemma \ref{lem:cj-d}}

\

\smallskip

\noindent
\textit{Case $d=1$.}
Invoking the definition \eqref{e:G}, one has here
\begin{align}\label{diff-g-1}
G_t(y-2^{-j}z)-G_t(y)= \frac{1}{2} \big\{\1_{\{|y-2^{-j}z|\leq t\}}-\1_{\{|y|\leq t\}}\big\} \,.
\end{align}
Whenever $|z|\le \min\{ 2^{j+1} t, 2\}$, a careful analysis of the intervals at stake reveals that 
\begin{eqnarray}\label{e:intervals}
\big| \1_{\{|y-2^{-j}z|\leq t\}}-\1_{\{|y|\leq t\}}\big|
&\lesssim& 
\1_{\{|y-t|\leq 2^{-j}|z|\}}+\1_{\{|y+t|\leq 2^{-j}|z|\}} \notag \\
&\lesssim& 
\1_{\{|y-t|\leq 2^{-j+1}\}}+\1_{\{|y+t|\leq 2^{-j+1}\}}\,. 
\end{eqnarray}
Furthermore,  if $2^{j+1}t<|z|<2$, then $t\leq 2^{-j}$ and
\begin{equation}\label{e:intervals'}
\big|\1_{\{|y-2^{-j}z|\leq t\}}-\1_{\{|y|\leq t\}}\big|\lesssim \1_{\{|y-2^{-j}z|\leq t\}}+\1_{\{|y|\leq t\}}\lesssim \1_{\{|y|\leq 2^{-j+2}\}} \, .
\end{equation}
Combining the bounds \eqref{e:intervals}-\eqref{e:intervals'} with the expression \eqref{diff-g-1}, we obtain
\begin{align*}
\big\|\ck^{(1)}_{t,j}\big\|_{L^1}&\leq 2^{j}\int_{\R}dy\int_{|z|<2} dz \, \big|G_t(y-2^{-j}z)-G_t(y)\big|\\
&\lesssim 2^{j}\int_{\R}dy\,  \big\{\1_{\{|y-t|\leq 2^{-j+1}\}}+\1_{\{|y+t|\leq 2^{-j+1}\}}+\1_{\{|y|\leq 2^{-j+2}\}}\big\}\\
&\lesssim 2^{j}\int_{\R}dy\,  \1_{\{|y|\leq 2^{-j+2}\}} \lesssim 2^{j}\int_{\R}dy\,  \1_{\{|2^jy|\leq 4\}} \lesssim \int_{\R}dy\,  \1_{\{|y|\leq 4\}} ,
\end{align*}
which naturally corresponds to the desired uniform bound.

\noindent
\textit{Case $d=2$.} Let us  denote the difference within the integral as
\begin{align}
A_{t, j,z}(y)& := G_t(y-2^{-j}z)-G_t(y)\nonumber\\
&=\frac{1}{\left( t^2 - |y- 2^{-j} z|^2\right)^{1/2}} \1_{\{t-|y-2^{-j}z|>0\} } - \frac{1}{\left( t^2 - |y|^2\right)^{1/2}} \1_{\{t-|y|>0\} }\,.\label{expr-a-tj}
\end{align}
For every $z\in \bfb_2$, we can write
\begin{align}
\int_{\R^2}dy \, |A_{t, j,z}(y) |&=\int_{\{|y-2^{-j}z|\leq |y|\}}dy \, |A_{t, j,z}(y)|+\int_{\{|y|\leq |y-2^{-j}z|\}}dy \, |A_{t, j,z}(y) |\nonumber\\
&=\int_{\{|y-2^{-j}z|\leq |y|\}}dy \, |A_{t, j,z}(y)|+\int_{\{|y+2^{-j}z|\leq |y|\}}dy \, |A_{t, j,z}(y+2^{-j}z) |\nonumber\\
&=\ci_{t,j}(z)+\ci_{t,j}(-z),\label{ci-z-ci-z}
\end{align}
where we have set
$$\ci_{t,j}(z):=\int_{\{|y-2^{-j}z|\leq |y|\}}dy \, |A_{t, j,z}(y)|.$$
Let us split the latter expression into
\begin{align}
\ci_{t,j}(z)&=\int_{\{t\leq |y-2^{-j}z|\leq |y|\}}dy \, |A_{t, j,z}(y)|\nonumber\\
&\hspace{1cm}+\int_{\{|y-2^{-j}z|<t\leq |y|\}}dy \, |A_{t, j,z}(y)|+\int_{\{|y-2^{-j}z|\leq |y|<t\}}dy \, |A_{t, j,z}(y)|\nonumber\\
&=:\ci^{(1)}_{t,j}(z)+\ci^{(2)}_{t,j}(z)+\ci^{(3)}_{t,j}(z).\label{decomp-ci}
\end{align}
We will now treat those three terms separately.

Given the expression \eqref{expr-a-tj} of $A_{t, j,z}(y)$, it is readily checked that $\ci^{(1)}_{t,j}(z)=0$. 
Next for the term $\ci^{(2)}_{t,j}(z)$, we notice that for every $z\in \bfb_2$,
\begin{align}
\ci^{(2)}_{t,j}(z)&= \int_{\{|y-2^{-j}z|<t\leq |y|\}}\frac{dy}{\left( t^2 - |y- 2^{-j} z|^2\right)^{1/2}}= \int_{\{|y|<t\leq |y+2^{-j}z|\}}\frac{dy}{\left( t^2 - |y|^2\right)^{1/2}}\label{bou-i-2-1}\\
&\leq  \int_{\{\max(0,t-2^{-j+1})<|y|<t\}}\frac{dy}{\left( t^2 - |y|^2\right)^{1/2}}\nonumber\\
&\lesssim \1_{\{0\leq t\leq 2^{-j+1}\}}  \int_0^t\frac{r\, dr}{\left( t^2 - r^2\right)^{1/2}}+ \1_{\{t> 2^{-j+1}\}}  \int_{t-2^{-j+1}}^t\frac{r\, dr}{\left( t^2 - r^2\right)^{1/2}}\nonumber\\
&\lesssim \1_{\{0\leq t\leq 2^{-j+1}\}}  \int_0^{t^2}\frac{ds}{\left( t^2 - s\right)^{1/2}}+ \1_{\{t> 2^{-j+1}\}}  \int_{(t-2^{-j+1})^2}^{t^2}\frac{ ds}{\left( t^2 - s\right)^{1/2}}\nonumber\\
&\lesssim \1_{\{0\leq t\leq 2^{-j+1}\}}  t+ \1_{\{t> 2^{-j+1}\}} (t^2-(t-2^{-j+1})^2)^{1/2}\lesssim 2^{-j+1}+ (2^{-j+1}t)^{1/2},\nonumber
\end{align}
which immediately entails that
\begin{equation}\label{bou-i-2-2}
\sup_{z\in \bfb_2}\sup_{t\in [0,1]}\ci^{(2)}_{t,j}(z)\lesssim 2^{-\frac{j}{2}}.
\end{equation}
Finally, for every $z\in \bfb_2$,
\begin{align}
\ci^{(3)}_{t,j}(z)&= \int_{\{|y-2^{-j}z|< |y|<t\}}dy\,\bigg|\frac{1}{\left( t^2 - |y- 2^{-j} z|^2\right)^{1/2}} - \frac{1}{\left( t^2 - |y|^2\right)^{1/2}} \bigg|\nonumber\\
&= \int_{\{|y-2^{-j}z|< |y|<t\}}dy\, \bigg[\frac{1}{\left( t^2 - |y|^2\right)^{1/2}}-\frac{1}{\left( t^2 - |y- 2^{-j} z|^2\right)^{1/2}} \bigg]\nonumber\\
&= \int_{\{|y-2^{-j}z|< |y|<t\}}\frac{dy}{\left( t^2 - |y|^2\right)^{1/2}}-\int_{\{|y|< |y+2^{-j}z|<t\}}\frac{dy}{\left( t^2 - |y|^2\right)^{1/2}}\nonumber\\
&=\int_{\R^2}\frac{dy}{\left( t^2 - |y|^2\right)^{1/2}} \big[ \1_{\{|y-2^{-j}z|< |y|<t\}}-\1_{\{|y|< |y+2^{-j}z|<t\}}\big].\label{bou-i-3-1}
\end{align}
At this point, the key observation is that
\begin{equation}\label{diff-indic}
\1_{\{|y-2^{-j}z|< |y|<t\}}-\1_{\{|y|< |y+2^{-j}z|<t\}}\leq \1_{\{|y|<t\leq |y+2^{-j}z|\}}.
\end{equation}
Indeed, the difference $\1_{\{|y-2^{-j}z|< |y|<t\}}-\1_{\{|y|< |y+2^{-j}z|<t\}}$ is strictly positive if and only if
$$|y-2^{-j}z|< |y|<t  \quad \text{and} \quad \big( |y| \geq |y+2^{-j}z| \quad \text{or} \quad |y+2^{-j}z|\geq t\big)$$
It turns out that we cannot have simultaneously $|y-2^{-j}z|< |y|$ and $ |y| \geq|y+2^{-j}z|$, since it would imply that $-2^{-j+1}y\cdot z+2^{-2j}|z|^2<0$ and $2^{-j+1}y\cdot z+2^{-2j}|z|^2\leq 0$, leading to an immediate contradiction. Thus, the difference $\1_{\{|y-2^{-j}z|< |y|<t\}}-\1_{\{|y|< |y+2^{-j}z|<t\}}$ is strictly positive if and only if
$$|y-2^{-j}z|< |y|<t  \quad \text{and} \quad  |y+2^{-j}z|\geq t,$$
which immediately yields \eqref{diff-indic}.

\smallskip

Injecting \eqref{diff-indic} into \eqref{bou-i-3-1}, we obtain for every $z\in \bfb_2$,
\begin{align*}
\ci^{(3)}_{t,j}(z)&\leq \int_{\R^2}\frac{dy}{\left( t^2 - |y|^2\right)^{1/2}}\1_{\{|y|<t\leq |y+2^{-j}z|\}}.
\end{align*}
We are here in the same position as in \eqref{bou-i-2-1}, and so we can use the same estimates as above to derive that uniformly over $z\in \bfb_2$ and $t\in [0,1]$,
\begin{equation}\label{bou-i-3-2}
\ci^{(3)}_{t,j}(z)\lesssim 2^{-\frac{j}{2}}.
\end{equation}
Plugging \eqref{bou-i-2-2} and \eqref{bou-i-3-2} into \eqref{decomp-ci}, we see that
$$\sup_{z\in \bfb_2}\sup_{t\in [0,1]} \ci_{t,j}(z) \lesssim 2^{-\frac{j}{2}},$$
which, going back to decomposition \eqref{ci-z-ci-z}, easily leads us to the assertion \eqref{asser-cj}.
\qed

\bigskip

\paragraph{\textbf{Acknowledgment.}} X. Chen is partially supported by the Simons Foundation grant 585506. J. Song is partially supported by Shandong University grant 111400\-89963041 and National Natural Science Foundation of China grant 12071256. S. Tindel is partially supported by the NSF grant  DMS-1952966.


\bigskip





\begin{thebibliography}{99}


\bibitem{balan22} R. M. Balan. Stratonovich Solution for the Wave Equation. {\it Journal of Theoretical Probability} (2022).

\bibitem{cdot} X. Chen, A, Deya, C, Ouyang, and S. Tindel. A $K$-rough path above the space-time fractional Brownian motion. {\it Stochastics and Partial Differential Equations: Analysis and Computations} (2021) Vol.31, 819-866.


\bibitem{cdot'} X. Chen, A, Deya, C, Ouyang, and S. Tindel. Moment estimates for some renormalized parabolic Anderson models. {\it The Annals of Probability} (2021) Vol.49, 2599-2636.



\bibitem{CDST}
X. Chen, A. Deya, J. Song and S. Tindel:
Solving the hyperbolic Anderson model 1: Skorohod setting

\bibitem{mourrat-weber}
J.-C. Mourrat and H. Weber: Global well-posedness of the dynamic $\Phi^4$ model in the plane. {\it Ann. Probab.} {\bf 45} (2017), no. 4.

\bibitem{deya19} A. Deya. A nonlinear wave equation with fractional perturbation. {\it Ann. Probab.} {\bf 47} (2019), no. 3, 1775-1810. 

\bibitem{gv95} J. Ginibre and G. Velo. Generalized Strichartz inequalities for the wave equation. {\it Journal of Functional Analysis} {\bf 133} (1995), 50-68. 

\bibitem{gko18} M. Gubinelli, H. Koch, and T. Oh. Renormalization of the two-dimensional stochastic nonlinear wave equations. {\it Trans. Amer. Math. Soc.} 370 (2018), 7335-7359.


\bibitem{rychkov}
V. S. Rychkov: Littlewood-Paley theory and function spaces with $A_p^{\textnormal{loc}}$ weights. {\it Math. Nachr.}, {\bf 224} (2001), no. 1.

\bibitem{qt07} L. Quer-Sardanyons and S. Tindel:   The 1-d stochastic wave equation driven by a fractional Brownian sheet. {\it Stochastic Processes and their Applications} {\bf 117} (2007), 1448-1472. 


\bibitem{st94}
G. Samorodnitsky and M.S. Taqqu: Stable non-gaussian random processes. Chapman and Hall, Boca Raton (1994)

\bibitem{stein93} E.M. Stein: Harmonic Analysis: Real – Variable Methods, Orthogonality, and Oscillatory Integrals, Princeton Univ. Press, Princeton, NJ, 1993


\end{thebibliography}
\end{document}